\documentclass{article}
\usepackage[utf8]{inputenc}
\usepackage{amssymb}
\usepackage{latexsym}
\usepackage{amsmath}
\usepackage{array}
\usepackage{graphicx}
\graphicspath{Images/}
\usepackage{caption}
\usepackage{subcaption}
\usepackage{algorithmic}
\usepackage{float}
\usepackage{hyperref}
\usepackage{cleveref}
\usepackage{algorithm}
\usepackage{amsthm}
\usepackage{verbatim}
\usepackage{color}
\usepackage{geometry}
\newtheorem{theorem}{Theorem}
\newtheorem{definition}{Definition}

\newtheorem{remark}{Remark}

\begin{document}

\title{On Lattice Boltzmann Methods based on vector-kinetic models for hyperbolic partial differential equations}
\author{Megala Anandan\footnote{Indian Institute of Science, C V Raman Road, Bangalore-560012, India. e-mail: megalaa@iisc.ac.in}, S. V. Raghurama Rao\footnote{Indian Institute of Science, C V Raman Road, Bangalore-560012, India. e-mail: raghu@iisc.ac.in}}
\date{January 8, 2024}

\maketitle

\begin{abstract}
In this paper, we are concerned about the lattice Boltzmann methods (LBMs) based on vector-kinetic models for hyperbolic partial differential equations. In addition to usual lattice Boltzmann equation (LBE) derived by explicit discretisation of vector-kinetic equation (VKE), we also consider LBE derived by semi-implicit discretisation of VKE and compare the relaxation factors of both. We study the properties such as H-inequality, total variation boundedness and positivity of both the LBEs, and infer that the LBE due to semi-implicit discretisation naturally satisfies all the properties while the LBE due to explicit discretisation requires more restrictive condition on relaxation factor compared to the usual condition obtained from Chapman-Enskog expansion. We also derive the macroscopic finite difference form of the LBEs, and utilise it to establish the consistency of LBEs with the hyperbolic system. Further, we extend this LBM framework to hyperbolic conservation laws with source terms, such that there is no spurious numerical convection due to imbalance between convection and source terms. We also present a D$2$Q$9$ model that allows upwinding even along diagonal directions in addition to the usual upwinding along coordinate directions. The different aspects of the results are validated numerically on standard benchmark problems.        
\end{abstract}

\textbf{Keywords: }
Vector-kinetic model, lattice Boltzmann method (LBM), H-inequality, total variation boundedness, positivity, consistency, source term, spurious numerical convection


\section{Introduction}
Lattice Boltzmann methods (LBMs) have emerged as a powerful and versatile class of computational techniques for simulating fluid flow and related phenomena. Over the years, they have gained significant popularity due to their ability to handle a wide range of fluid flow scenarios, from incompressible flows (\cite{Incomp1,Incomp2,Incomp3}) to complex multiphase (\cite{Multiphase1,Multiphase2,Multiphase3,Multiphase4,Multiphase5}) and multiscale (\cite{Multiscale1}) systems. LBMs have been employed for modelling and simulating problems in magnetohydrodynamics (\cite{MHD1,MHD2,MHD3,MHD4}), porous media (\cite{Por1,Por2,Por3,Por4}), heat transfer (\cite{Heattransfer1,Heattransfer2,Heattransfer3}) and turbulence (\cite{Turb1,Turb2}). The reader is referred to the books \cite{Basic_mohamad,Basic_succi,Basic_Guo_Shu} for extensive study of LBMs, \cite{Basic_chen} for review of LBMs for fluid flows, \cite{Basic_DAPerumal_AnoopDas} for review of LBMs for heat transfer, and \cite{Basic_Ansumali} for review of entropic LBMs. \\
The Lattice Boltzmann equation (LBE) has been shown to approximate the Euler and the Navier-Stokes equations through different approaches such as Chapman-Enskog expansion (\cite{CEanal3,CEanal4,CEanal5}), asymptotic expansion (\cite{AEanal1,AEanal2,AEanal3}), Maxwellian iteration (\cite{Maxiter1,Maxiter2,Maxiter3}), equivalent equation (\cite{EqeqnFDubois}), and recursive representation (\cite{Recrep1}). Various attempts have been made in which the LBE is shown to be equivalent to mutistep finite difference equation (\cite{Suga_FD,Dellacherie_FD,Ginzburg_FD,Fucik_FD,Graille_FD}), and the consistency with macroscopic equations has been shown in \cite{Bellotti_FDcons}. Further, the linkage between LBM and relaxation systems of \cite{RS_Jin_Xin} was explored in \cite{RS_DV_Rheinlander,RS_DV_Graille,RS_DV_Krause}.  \\
While the discussions above correspond to the LBE derived from discretisation of the Boltzmann equation (essentially scalar-kinetic equation) with discrete velocities, we consider the class of LBEs derived from discretisation of vector-kinetic equations introduced in \cite{Bouchut1999,Bouchut2003,10.2307/2587356}. The vector-kinetic models have been utilised to develop various numerical schemes in the areas of porous media \cite{VK_Jobic}, entropy stable methods for hyperbolic systems \cite{VK_Megala}, implicit kinetic relaxation schemes \cite{VK_Courtes}, and lattice Boltzmann relaxation schemes \cite{VK_Rao,VK_Rohan_thesis,VK_Ankit}. In particular, \cite{VK_Courtes} and \cite{VK_Rao} present the lattice Boltzmann methods with different equilibrium functions and their resulting Chapman-Enskog expansions. In this paper, we present some important properties (such as macroscopic multi-step finite difference form and consistency) of the LBE derived from vector-kinetic equations. We also present a novel way to handle well-balancing of convection and source terms in this framework. Further, we also present an LBM model that allows upwinding along diagonal directions in addition to the usual upwinding along coordinate directions (presented first in the proceedings of a conference \cite{D2Q9_Megala}).  \\ 
The paper is organised as follows: Section \ref{CH_LB_SecMM} presents the mathematical model of hyperbolic conservation law and its vector-kinetic equation. Section \ref{Ch_LB_SecLB} presents two different ways of deriving LBE from vector-kinetic equation, their Chapman-Enskog expansion and different equilibrium functions. Different properties such as H-inequality, macroscopic multi-step finite difference form, consistency, total variation boundedness and positivity are discussed in section \ref{Ch_LB_SecProp}. The well-balancing technique for hyperbolic partial differential equations with source terms is explained in section \ref{Ch_LB_SecSource}. The D$2$Q$9$ model of LBM that allows upwinding along diagonal directions is explained in section \ref{Ch_LB_SecD2Q9}. The numerical validation of the methods is presented in section \ref{Ch_LB_SecNum}. Section \ref{Ch_LB_SecConc} concludes the paper. 

\section{Mathematical model}
\label{CH_LB_SecMM}
In this section, we describe the hyperbolic conservation law and the vector-kinetic equation that approximates it.
\subsection{Hyperbolic conservation law}
Consider the hyperbolic conservation law
\begin{equation}
\label{Ch_LB_Hyp cons law}
    \partial_t U + \partial_{x_d} G^d(U)=0,
\end{equation}
where $U(x_1,x_2,\dots,x_D,t):\Omega \times [0,T] \to \mathbb{R}^p$ is the conserved variable and $G^d(U):\mathbb{R}^p \to \mathbb{R}^p$ is the flux in direction $d$, for $d \in \{1,2,\dots,D \}$. Here $\Omega \subset \mathbb{R}^D$, $D$ and $p$ indicate the number of dimensions and number of equations in the system respectively. $\eta(U)$ is the convex entropy function for \eqref{Ch_LB_Hyp cons law}. 
\subsection{Vector-kinetic equation}
The hyperbolic conservation law in \eqref{Ch_LB_Hyp cons law} can be approximated by the vector-kinetic equation (VKE), 
\begin{equation}
\label{Ch_LB_DVBE}
    \partial_t f_q + \partial_{x_d} \left( v^d_q f_q \right)= - \frac{1}{\epsilon} \left( f_q - f^{eq}_q(U) \right).
\end{equation}
Here $f_q:\Omega \times [0,T] \to \mathbb{R}^p$, $f^{eq}_q:\mathbb{R}^p \to \mathbb{R}^p$ and $q \in \{1,2,\dots,Q \}$ with $Q$ being the number of discrete velocities. $\epsilon$ is a positive small parameter. $v^d_q$ is the $d^{th}$ component of the $q^{th}$ discrete velocity. Summing \eqref{Ch_LB_DVBE} over all $q$, we get
\begin{equation}
    \partial_t \sum_{q=1}^Q f_q + \partial_{x_d} \sum_{q=1}^Q \left( v^d_q f_q \right)= - \frac{1}{\epsilon} \sum_{q=1}^Q \left( f_q - f^{eq}_q(U) \right).
\end{equation}
If $\sum_{q=1}^Q f_q = \sum_{q=1}^Q f^{eq}_q = U$, then
\begin{equation}
\label{Ch_LB_DVBE_sum}
    \partial_t U + \partial_{x_d} \sum_{q=1}^Q \left( v^d_q f_q \right)= 0.
\end{equation}
In the limit $\epsilon \to 0$, we infer from \eqref{Ch_LB_DVBE} that $f_q \to f^{eq}_q(U)$. Thus, we can write $f_q$ as perturbation (in $\epsilon$) of $f^{eq}_q$:
\begin{equation}
\label{Ch_LB_DVBE_neq}
    f_q=f^{eq}_q + \epsilon f^{neq}_q,
\end{equation}
where $f^{neq}_q$ consists of the non-equilibrium perturbations. \\
If $\sum_{q=1}^Q v^d_q f^{eq}_q = G^d(U)$, then \eqref{Ch_LB_DVBE_sum} becomes the hyperbolic conservation law \eqref{Ch_LB_Hyp cons law} in the limit $\epsilon \to 0$. 

\section{Lattice Boltzmann equation}
\label{Ch_LB_SecLB}
In this section, we present explicit and semi-implicit lattice Boltzmann discretisations of the VKE \eqref{Ch_LB_DVBE}, their comparison, and their Chapman-Enskog expansions. \\
Let us use the vector notations: $\mathbf{x}=\begin{bmatrix}
    x_1, & x_2, & \dots, & x_D \end{bmatrix}$ and $\mathbf{v_q}=\begin{bmatrix}
    v^1_q, & v^2_q, & \dots, & v^D_q \end{bmatrix}$. An explicit Euler discretisation of the VKE \eqref{Ch_LB_DVBE} along $\frac{d x_d}{dt}=v^d_q$ (the characteristic equation) gives 
\begin{equation}
    f_q\left( \mathbf{x},t+\Delta t \right) = f_q \left( \mathbf{x}-\mathbf{v_q}\Delta t,t\right) - \frac{\Delta t}{\epsilon} \left( f_q \left( \mathbf{x}-\mathbf{v_q}\Delta t,t\right) - f^{eq}_q \left( U \left( \mathbf{x}-\mathbf{v_q}\Delta t,t\right) \right)  \right).
\end{equation}
Using $\omega=\frac{\Delta t}{\epsilon}$ and rewriting the above equation, we obtain the lattice Boltzmann equation (LBE)
\begin{equation}
\label{Ch_LB_LBEexplicit}
    f_q\left( \mathbf{x},t+\Delta t \right) = (1-\omega) f_q \left( \mathbf{x}-\mathbf{v_q}\Delta t,t\right) +\omega f^{eq}_q \left( U \left( \mathbf{x}-\mathbf{v_q}\Delta t,t\right) \right).  
\end{equation}
On the other hand, a semi-implicit discretisation of the VKE \eqref{Ch_LB_DVBE} with  implicit treatment of $f_q$ in the collision term gives
\begin{equation}
\label{Ch_LB_semiimplicitdisc}
    f_q\left( \mathbf{x},t+\Delta t \right) = f_q \left( \mathbf{x}-\mathbf{v_q}\Delta t,t\right) - \frac{\Delta t}{\epsilon} \left( f_q \left( \mathbf{x},t+\Delta t \right) - f^{eq}_q \left( U \left( \mathbf{x}-\mathbf{v_q}\Delta t,t\right) \right)  \right).
\end{equation}   
Rewriting the above equation as
\begin{gather}
    f_q\left( \mathbf{x},t+\Delta t \right) = \left( \frac{1}{1+\omega} \right) f_q \left( \mathbf{x}-\mathbf{v_q}\Delta t,t\right) + \left( \frac{\omega}{1+\omega} \right) f^{eq}_q \left( U \left( \mathbf{x}-\mathbf{v_q}\Delta t,t\right) \right), \\
    \label{Ch_LB_LBEsemiimplicit}
   \textrm{or} \ f_q\left( \mathbf{x},t+\Delta t \right) = \left(1-\Tilde{\omega}\right) f_q \left( \mathbf{x}-\mathbf{v_q}\Delta t,t\right) +\Tilde{\omega} f^{eq}_q \left( U \left( \mathbf{x}-\mathbf{v_q}\Delta t,t\right) \right),
\end{gather}
an LBE with $\Tilde{\omega}=\frac{\omega}{1+\omega}$ is obtained. \\
If the grid is uniform with spacing $\Delta x_d$ along direction $d$ and if the velocities are chosen such that $v^d_q=m\frac{\Delta x_d}{\Delta t}$ with $m \in \mathbb{Z}$, $\forall d,q$, then the collision-streaming algorithm
\begin{eqnarray}
    \label{Ch_LB_Collision}
    \text{Collision:} & f_q^*\left( \mathbf{x}-\mathbf{v_q}\Delta t,t\right) = (1-\hat{\omega}) f_q \left( \mathbf{x}-\mathbf{v_q}\Delta t,t\right) +\hat{\omega} f^{eq}_q \left( U \left( \mathbf{x}-\mathbf{v_q}\Delta t,t\right) \right) \\
    \label{Ch_LB_Streaming}
    \text{Streaming:} & f_q\left( \mathbf{x},t+\Delta t \right) = f_q^*\left( \mathbf{x}-\mathbf{v_q}\Delta t,t\right) 
\end{eqnarray}
can be used to numerically implement the LBEs in \eqref{Ch_LB_LBEexplicit} (with $\hat{\omega}=\omega$) and \eqref{Ch_LB_LBEsemiimplicit} (with $\hat{\omega}=\Tilde{\omega}$). It is to be noted that the streaming in \eqref{Ch_LB_Streaming} is exact. After evaluating $f_q\left( \mathbf{x},t+\Delta t \right)$, we find $U$ by using $U\left( \mathbf{x},t+\Delta t \right) = \sum_{q=1}^Q f_q\left( \mathbf{x},t+\Delta t \right)$. Then, we evaluate $f^{eq}_q \left( U \left( \mathbf{x},t+\Delta t\right) \right)$ and then proceed with the next time step. 
{\em Hereafter, we use $\hat{\omega}$ in the presentation of our theory to commonly represent $\omega$ in \eqref{Ch_LB_LBEexplicit} and $\Tilde{\omega}$ in \eqref{Ch_LB_LBEsemiimplicit}}.   

\subsection{Chapman-Enskog expansion}
Taylor expanding the LBEs in \eqref{Ch_LB_LBEexplicit} (with $\hat{\omega}=\omega$) and \eqref{Ch_LB_LBEsemiimplicit} (with $\hat{\omega}=\Tilde{\omega}$) and simplifying, we get
\begin{equation}
\label{Ch_LB_CE_TE2}
\left( \partial_t + \mathbf{v_q} \cdot \nabla \right) f_q =  - \frac{\hat{\omega}}{\Delta t} \left(f_q-f^{eq}_q\right) + \frac{\hat{\omega}}{2} \left( \partial_t +  \mathbf{v_q} \cdot \nabla \right) \left(f_q-f_q^{eq}\right) + \mathcal{O}(\Delta t^2).
\end{equation}
Consider the perturbation expansion of $f_q$:
\begin{equation}
\label{Ch_LB_CE_PE}
f_q=f_q^{eq} + \epsilon f_q^{(1)} + \epsilon^2 f_q^{(2)} + \dots
\end{equation}
Using the above expression, since $\sum_{q=1}^Q f_q = \sum_{q=1}^Q f_q^{eq} = U$, we infer that the moment of non-equilibrium distribution function leads to $\sum_{q=1}^Q \left(\epsilon f_q^{(1)} + \epsilon^2 f_q^{(2)}+ \dots \right)=0$. Each term corresponding to different order of $\epsilon$ in this moment expression must individually be zero. Hence $\sum_{q=1}^Q f_q^{(i)} = 0, \forall i \in \mathbb{N}$. Multiple scale expansion of derivatives of $f_q$ gives $\partial_t f_q = \left( \epsilon \partial_t^{(1)} + \epsilon^2 \partial_t^{(2)} + ... \right) f_q \ \text{and} \ \mathbf{v_q} \cdot \nabla f_q = \epsilon \mathbf{v_q} \cdot \nabla^{(1)} f_q$. \\
Using perturbation expansion of $f_q$ and multiple scale expansion of derivatives of $f_q$ in \eqref{Ch_LB_CE_TE2} and separating out $O(\epsilon)$ and $O(\epsilon^2)$ terms,
 \begin{eqnarray}
 \label{Ch_LB_CE_O(xi)}
 O(\epsilon): & \left( \partial_t^{(1)} + \mathbf{v_q} \cdot \nabla^{(1)}\right) f_q^{eq} = - \frac{\hat{\omega}}{\Delta t} f_q^{(1)} \\
 \label{Ch_LB_CE_O(xi^2)}
 O(\epsilon^2): & \ \partial_t^{(2)} f_q^{eq} + \left( 1- \frac{\hat{\omega}}{2}\right)\left( \partial_t^{(1)} + \mathbf{v_q} \cdot \nabla^{(1)}\right) f_q^{(1)} =  - \frac{\hat{\omega}}{\Delta t} f_q^{(2)}  
 \end{eqnarray}
Zeroth moment $\left(\sum_{q=1}^Q\right)$ of $O(\epsilon)$ terms in \eqref{Ch_LB_CE_O(xi)} and $O(\epsilon^2)$ terms in \eqref{Ch_LB_CE_O(xi^2)} respectively give
 \begin{gather}
 \label{Ch_LB_CE_0thO(xi)}
 \partial_t^{(1)} U + \partial_{x_d}^{(1)} G^d(U) = 0,  \\
 \label{Ch_LB_CE_0thO(xi^2)}
 \partial_t^{(2)} U + \left( 1- \frac{\hat{\omega}}{2}\right) \partial_{x_d}^{(1)} \left(\sum_{q=1}^Q v_q^{(d)}  f_q^{(1)}\right)=0.  
 \end{gather}
From the first moment $\left(\sum_{q=1}^Q v^d_q\right)$ of $O(\epsilon)$ terms in \eqref{Ch_LB_CE_O(xi)}, we get 
\begin{equation}
\label{Ch_LB_CE_1stO(xi)}
\sum_{q=1}^Q v_q^d f_q^{(1)} = - \frac{\Delta t}{\hat{\omega}} \left( \partial_U G^d \left( -\partial_U G^i \partial_{x_i}^{(1)} U \right) + \partial_{x_i}^{(1)} \left( \sum_{q=1}^Q v_q^d v_q^i f_q^{eq} \right)\right)
\end{equation}
Recombining the zeroth moment equations of $O(\epsilon)$ in \eqref{Ch_LB_CE_0thO(xi)} and $O(\epsilon^2)$ in \eqref{Ch_LB_CE_0thO(xi^2)} and reversing the multiple scale expansions, we get 
\begin{equation}
\label{Ch_LB_CE_mPDE2}
\partial_t U + \partial_{x_d} G^d(U) = \\ \Delta t \left( \frac{1}{\hat{\omega}}-\frac{1}{2}\right) \partial_{x_d} \left(  \partial_{x_i} \left( \sum_{q=1}^Q v_q^d v_q^i f_q^{eq} \right) - \partial_U G^d \partial_U G^i \partial_{x_i} U \right)
\end{equation}
upto $\mathcal{O}\left( \Delta t^2 \right)$. 
\subsection{Equilibrium function}
\label{Ch_LB_Sec_Eqbm fun}
In the previous sections, we imposed the following conditions on $f^{eq}_q$:
\begin{equation}
\label{Ch_LB_eqb_moms}
    \sum_{q=1}^Q f^{eq}_q = U, \ \sum_{q=1}^Q v^d_q f^{eq}_q = G^d(U).
\end{equation}
In this section, we present some $f^{eq}_q$ that satisfy the above requirements. 
\subsubsection{Classical D$1$Q$2$}
Consider one dimension (D=1) and 2 discrete velocities ($Q=2$) such that $v^1_1=\lambda$ and $v^1_2=-\lambda$, and $\lambda=\frac{\Delta x_1}{\Delta t}$. Then, 
\begin{equation}
    f^{eq}_q=\frac{1}{2}U - (-1)^q \frac{1}{2\lambda} G^1(U) \text{ for } q \in \{ 1,2\}
\end{equation}
satisfies \eqref{Ch_LB_eqb_moms}. The Chapman-Enskog expansion \eqref{Ch_LB_CE_mPDE2} in this case becomes,
\begin{equation}
    \partial_t U + \partial_{x_1} G^1(U) = \Delta t \left( \frac{1}{\hat{\omega}}-\frac{1}{2}\right) \partial_{x_1} \left( \left( \lambda^2 I - |\partial_U G^1|^2 \right) \partial_{x_1} U \right). 
\end{equation}
It is to be noted that the $\mathcal{O}(\Delta t)$ term on the right hand side of the above equation represents numerical diffusion. For stability, we require the numerical diffusion coefficient to be positive. Therefore, we require $\lambda^2 I > |\partial_U G^1|^2$ and $0<\hat{\omega}<2$. 
\subsubsection{D$1$Q$3$}
Consider one dimension (D=1) and 3 discrete velocities ($Q=3$) such that $v^1_1=\lambda$, $v^1_2=0$ and $v^1_3=-\lambda$, and $\lambda=\frac{\Delta x_1}{\Delta t}$. Then,
\begin{equation}
    f^{eq}_q=\frac{1}{3}U + \left( \delta_{q1} - \delta_{q3} \right) \frac{1}{2\lambda} G^1(U) \text{ for } q \in \{ 1,2,3\}
\end{equation}
where $\delta$ is the \textit{Kronecker delta function}, satisfies \eqref{Ch_LB_eqb_moms}. The Chapman-Enskog expansion \eqref{Ch_LB_CE_mPDE2} in this case becomes,
\begin{equation}
    \partial_t U + \partial_{x_1} G^1(U) = \\ \Delta t \left( \frac{1}{\hat{\omega}}-\frac{1}{2}\right) \partial_{x_1} \left( \left( \frac{2}{3} \lambda^2 I - |\partial_U G^1|^2 \right) \partial_{x_1} U \right). 
\end{equation}
Enforcement of the positivity of numerical diffusion coefficient yields $\lambda^2 I > \frac{3}{2}|\partial_U G^1|^2$ and $0<\hat{\omega}<2$.
\subsubsection{Upwind D$\overline{d}$Q$(2\overline{d}+1)$}
Consider $D=\overline{d}$ and $Q=2\overline{d}+1$ with $\lambda_d=\frac{\Delta x_d}{\Delta t}$ and
\begin{equation}
    v^d_q= \lambda_d \delta_{qd} - \lambda_d \delta_{q\left(d+\left(\overline{d}+1\right)\right)}.
\end{equation}
Define
\begin{eqnarray}
\label{Ch_LB_upweqbmfun}
    f^{eq}_q &=& \left\{ \begin{matrix} \frac{G^{q+}}{\lambda_q}, & \text{for } q \in \{1,2,\dots,\overline{d} \} \\
    U - \sum_{d=1}^{\overline{d}} \left( \frac{G^{d+}+G^{d-}}{\lambda_d} \right),  & \text{for } q=\overline{d}+1 \\
    \frac{G^{\left(q-\left(\overline{d}+1\right)\right)-}}{\lambda_{q-\left(\overline{d}+1\right)}}, & \text{for } q \in \{\overline{d}+2,\overline{d}+3,\dots,2\overline{d}+1 \}  \end{matrix} \right.
\end{eqnarray}
with $G^d=G^{d+}-G^{d-}$. This satisfies \eqref{Ch_LB_eqb_moms} and leads to the {\em Flux Decomposition} technique of \cite{10.2307/2587356}. Using an additional choice, $G^{d+}$ and $G^{d-}$ for a hyperbolic system can be evaluated by a suitable flux splitting method available in literature. For instance, one can use $G^{d+}$ and $G^{d-}$ from commonly known flux vector splitting methods such as kinetic flux vector splitting \cite{FVS_KFVS}, Steger-Warming flux vector splitting \cite{FVS_SW} and van Leer's flux vector splitting \cite{FVS_vL}. One can also evaluate $G^{d+}$ and $G^{d-}$ from some flux difference splitting methods such as Roe's approximate Riemann solver \cite{FVS_Roe} and kinetic flux difference splitting \cite{SHRINATH2023105702}. If we consider scalar conservation laws (\textit{i.e.,} $p=1$), then we can simply use the sign of wave speed $\partial_U G^d$ to determine the split fluxes as: 
\begin{eqnarray}
    \partial_UG^{d+} = \left\{ \begin{matrix}
        \partial_UG^d & \text{if } \partial_U G^d > 0 \\ 0 & \text{if } \partial_U G^d \leq 0
    \end{matrix}   \right., & \partial_UG^{d-} = \left\{ \begin{matrix}
        0 & \text{if } \partial_U G^d > 0 \\ -\partial_UG^d & \text{if } \partial_U G^d \leq 0
    \end{matrix}   \right., 
\end{eqnarray}
\begin{equation}
    G^{d\pm} = \int_{0}^U \partial_UG^{d\pm} dU \text{ if } G^d(U=0)=0. 
\end{equation}
The Chapman-Enskog expansion \eqref{Ch_LB_CE_mPDE2} for the case of upwind D$\overline{d}$Q$(2\overline{d}+1)$ becomes,
\begin{equation}
\label{Ch_LB_Eqbmupwind_CEexp}
    \partial_t U + \partial_{x_d} G^d(U) =  \Delta t \left( \frac{1}{\hat{\omega}}-\frac{1}{2}\right) \partial_{x_d} \left( \delta_{di} \lambda_d\partial_U \left( G^{d+} + G^{d-} \right) - \partial_U G^d \partial_U G^i \right)\partial_{x_i} U.
\end{equation}
For positivity of numerical diffusion coefficient, we require 
\begin{equation}
\label{Ch_LB_Upwindeqbm_subcharcondn}
    \lambda_d\partial_U \left( G^{d+} + G^{d-} \right) - \partial_U G^d \partial_U G^i > 0
\end{equation} 
along with $0<\hat{\omega}<2$. \\ 
For all the models of equilibrium function described above, a condition relating $\lambda_d$ and $\partial_U G^d$ is obtained while ensuring positivity of numerical diffusion coefficient. Such relations are known as sub-characteristic conditions as they relate the characteristic speeds of vector-kinetic equation to those of the hyperbolic conservation law.
\begin{remark}
In all the models of equilibrium function described above, $0<\hat{\omega}<2$ is required for enforcing the positivity of numerical diffusion coefficient. We know that $\hat{\omega}=\omega$ and $\hat{\omega}=\Tilde{\omega}$ for LBEs in \eqref{Ch_LB_LBEexplicit} and \eqref{Ch_LB_LBEsemiimplicit} respectively. Thus, the stability requirement is,
\begin{eqnarray}
    \text{For LBE in } \eqref{Ch_LB_LBEexplicit}: & 0<\hat{\omega}=\omega=\frac{\Delta t}{\epsilon}<2 \implies 0<\Delta t<2\epsilon, \\
    \text{For LBE in } \eqref{Ch_LB_LBEsemiimplicit}: & 0<\hat{\omega}=\Tilde{\omega}=\frac{\omega}{1+\omega}=\frac{\Delta t}{\epsilon+\Delta t}<2 \implies \Delta t > 0. 
\end{eqnarray}
It is to be noted that the requirement of $0<\Tilde{\omega}<2$ does not impose any upper-bound on $\Delta t$ for the LBE in \eqref{Ch_LB_LBEsemiimplicit}. 
\end{remark}

\begin{figure}[t]
\centering
\begin{subfigure}[b]{0.6\textwidth}
\centering
\includegraphics[width=\textwidth]{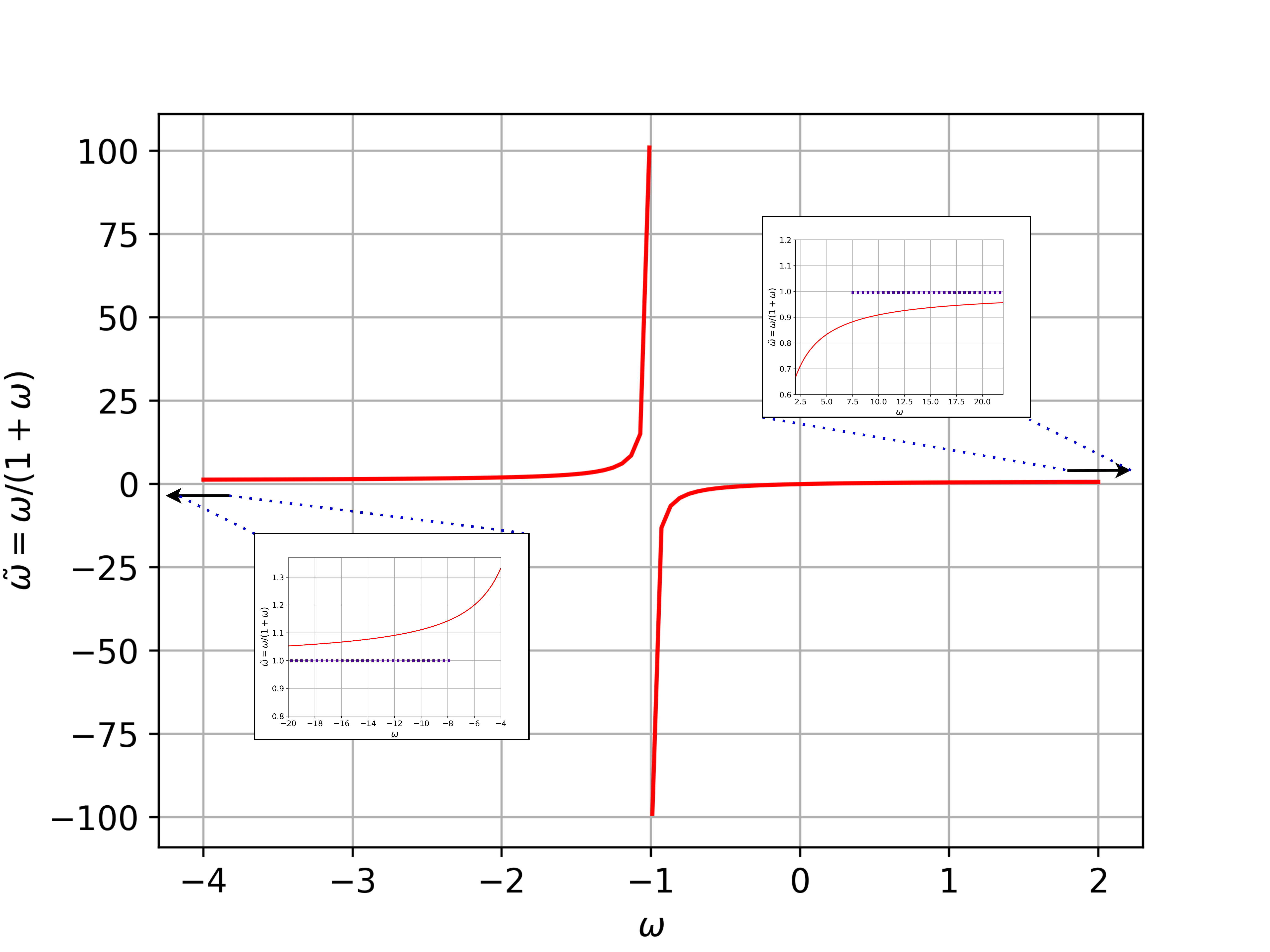}
\end{subfigure}
\caption{\centering Plot of $\Tilde{\omega}$ vs.  $\omega$}
\label{Fig:Omega_graph}
\end{figure}

\begin{remark}
    For the LBE in \eqref{Ch_LB_LBEsemiimplicit}, the positivity of numerical diffusion coefficient enforces $0<\Tilde{\omega}<2$. However, since $\Tilde{\omega}=\frac{\omega}{1+\omega}$ and $\omega=\frac{\Delta t}{\epsilon}>0$, $\Tilde{\omega}$ is restricted to the interval $(0,1)$. Figure \ref{Fig:Omega_graph} shows the plot of $\Tilde{\omega}$ vs.  $\omega$, and it can be seen that $0 < \Tilde{\omega} < 1$ for $\omega > 0$. 
\end{remark}

\section{Properties of the lattice Boltzmann equation}
\label{Ch_LB_SecProp}
In this section, we discuss the properties of LBEs in \eqref{Ch_LB_LBEexplicit} (with $\hat{\omega}=\omega$) and \eqref{Ch_LB_LBEsemiimplicit} (with $\hat{\omega}=\Tilde{\omega}$). The properties considered are: H-inequality, macroscopic finite difference form, consistency, total variation boundedness and positivity. 
\subsection{H-inequality}
We prove that an H-inequality is associated with the LBE obtained from semi-implicit discretisation of the VKE (\textit{i.e.,} \eqref{Ch_LB_LBEsemiimplicit}). We also show that a constraint on $\omega$ is required to associate an H-inequality with the LBE obtained from explicit discretisation of the VKE (\textit{i.e.,} \eqref{Ch_LB_LBEexplicit}). For convenience, we consider scalar conservation laws (\textit{i.e.,} $p=1$) in the presentation of H-inequality.  
\begin{definition}
\label{Ch_LB_Hineq_def}
    Define a function $H_q(f_q)$ such that:
    \begin{itemize}
        \item $H_q(f_q)$ is convex with respect to $f_q$ (\textit{i.e.,} $\frac{\partial H_q}{\partial f_q}$ is monotonically increasing and $\frac{\partial^2 H_q}{\partial f^2_q}$ is positive-definite),
        \item $\sum_{q=1}^Q H_q(f^{eq}_q) = \eta (U)$,
        \item $\sum_{q=1}^Q H_q(f^{eq}_q) \leq \sum_{q=1}^Q H_q(f_q)$. 
    \end{itemize}
\end{definition}
We consider the semi-implicit discretisation \eqref{Ch_LB_semiimplicitdisc} of VKE  with the notation $f_{q_{x_i}}^{n+1}:=f_q\left( \mathbf{x},t+\Delta t \right)$, $f_{q_{y_i}}^{n}:=f_q \left( \mathbf{x}-\mathbf{v_q}\Delta t,t\right)$ and $f_{q_{y_i}}^{eq^n}:=f^{eq}_q \left( \mathbf{x}-\mathbf{v_q}\Delta t,t\right)$:
\begin{equation}
\label{Ch_LB_H_fexp}
    f_{q_{x_i}}^{n+1} = f_{q_{y_i}}^{n} - \omega \left( f_{q_{x_i}}^{n+1} - f_{q_{y_i}}^{eq^n} \right).
\end{equation}
\begin{theorem}
     There exists an inequality 
    \begin{equation}
        H_q \left( f_{q_{x_i}}^{n+1} \right) - H_q \left(f_{q_{y_i}}^{n} \right) \leq - \omega \left( H_q \left( f_{q_{x_i}}^{n+1} \right) - H_q \left(f_{q_{y_i}}^{eq^n} \right) \right)
    \end{equation}
    corresponding to the semi-implicit discretisation \eqref{Ch_LB_H_fexp} of VKE with $\omega=\frac{\Delta t}{\epsilon}> 0$. Here, $H_q(f_q)$ follows the definition \ref{Ch_LB_Hineq_def}. 
\end{theorem}
\begin{proof}
    Left multiplying $\left. \frac{\partial H_q}{\partial f_q} \right\rvert_{f_{q_{x_i}}^{n+1}}$ to \eqref{Ch_LB_H_fexp}, we get
    \begin{equation}
    \label{Ch_LB_Hineq_proof}
        \left. \frac{\partial H_q}{\partial f_q} \right\rvert_{f_{q_{x_i}}^{n+1}} \left( f_{q_{x_i}}^{n+1}- f_{q_{y_i}}^{n} \right) = -\omega \left. \frac{\partial H_q}{\partial f_q} \right\rvert_{f_{q_{x_i}}^{n+1}} \left( f_{q_{x_i}}^{n+1} - f_{q_{y_i}}^{eq^n} \right).
    \end{equation} 
    We consider the left and right hand sides of the above equation separately. \\
    By mean value theorem, we have  
    \begin{equation}
        \left. \frac{\partial H_q}{\partial f_q} \right\rvert_{f_a} \left( f_{q_{x_i}}^{n+1}- f_{q_{y_i}}^{n} \right) = H_q \left( f_{q_{x_i}}^{n+1} \right) - H_q \left(f_{q_{y_i}}^{n} \right) 
    \end{equation}
    for some $f_a$ lying on the line segment connecting $f_{q_{x_i}}^{n+1}$ and $f_{q_{y_i}}^{n}$. Further, we have the following due to the monotonicity of $\frac{\partial H_q}{\partial f_q}$:
    \begin{gather}
        f_{q_{x_i}}^{n+1} \geq f_a \geq f_{q_{y_i}}^{n} \implies \left. \frac{\partial H_q}{\partial f_q} \right\rvert_{f_{q_{x_i}}^{n+1}} \geq \left. \frac{\partial H_q}{\partial f_q} \right\rvert_{f_a} \geq \left. \frac{\partial H_q}{\partial f_q} \right\rvert_{f_{q_{y_i}}^{n}}, \\
        f_{q_{x_i}}^{n+1} \leq f_a \leq f_{q_{y_i}}^{n} \implies \left. \frac{\partial H_q}{\partial f_q} \right\rvert_{f_{q_{x_i}}^{n+1}} \leq \left. \frac{\partial H_q}{\partial f_q} \right\rvert_{f_a} \leq \left. \frac{\partial H_q}{\partial f_q} \right\rvert_{f_{q_{y_i}}^{n}}. 
    \end{gather}
    Thus, we obtain the following inequality involving the term on the left hand side of \eqref{Ch_LB_Hineq_proof}:  
    \begin{equation}
    \label{Ch_LB_Hineq_proof_a}
        H_q \left( f_{q_{x_i}}^{n+1} \right) - H_q \left(f_{q_{y_i}}^{n} \right) = \left. \frac{\partial H_q}{\partial f_q} \right\rvert_{f_a} \left( f_{q_{x_i}}^{n+1}- f_{q_{y_i}}^{n} \right) 
        \leq \left. \frac{\partial H_q}{\partial f_q} \right\rvert_{f_{q_{x_i}}^{n+1}} \left( f_{q_{x_i}}^{n+1}- f_{q_{y_i}}^{n} \right). 
    \end{equation}
    On the other hand, we also have the following by mean value theorem:
    \begin{equation}
        \left. \frac{\partial H_q}{\partial f_q} \right\rvert_{f_b} \left( f_{q_{x_i}}^{n+1}- f_{q_{y_i}}^{eq^n} \right) = H_q \left( f_{q_{x_i}}^{n+1} \right) - H_q \left(f_{q_{y_i}}^{eq^n} \right) 
    \end{equation}
    for some $f_b$ lying on the line segment connecting $f_{q_{x_i}}^{n+1}$ and $f_{q_{y_i}}^{eq^n}$. Further, due to the monotonicity of $\frac{\partial H_q}{\partial f_q}$, we have
    \begin{gather}
        f_{q_{x_i}}^{n+1} \geq f_b \geq f_{q_{y_i}}^{eq^n} \implies \left. \frac{\partial H_q}{\partial f_q} \right\rvert_{f_{q_{x_i}}^{n+1}} \geq \left. \frac{\partial H_q}{\partial f_q} \right\rvert_{f_b} \geq \left. \frac{\partial H_q}{\partial f_q} \right\rvert_{f_{q_{y_i}}^{eq^n}}, \\
        f_{q_{x_i}}^{n+1} \leq f_b \leq f_{q_{y_i}}^{eq^n} \implies \left. \frac{\partial H_q}{\partial f_q} \right\rvert_{f_{q_{x_i}}^{n+1}} \leq \left. \frac{\partial H_q}{\partial f_q} \right\rvert_{f_b} \leq \left. \frac{\partial H_q}{\partial f_q} \right\rvert_{f_{q_{y_i}}^{eq^n}}. 
    \end{gather}
    Thus, we obtain the following inequality involving the term on the right hand side of \eqref{Ch_LB_Hineq_proof}:  
    \begin{equation}
    \label{Ch_LB_Hineq_proof_b}
        H_q \left( f_{q_{x_i}}^{n+1} \right) - H_q \left(f_{q_{y_i}}^{eq^n} \right) = \left. \frac{\partial H_q}{\partial f_q} \right\rvert_{f_b} \left( f_{q_{x_i}}^{n+1}- f_{q_{y_i}}^{eq^n} \right) 
        \leq \left. \frac{\partial H_q}{\partial f_q} \right\rvert_{f_{q_{x_i}}^{n+1}} \left( f_{q_{x_i}}^{n+1}- f_{q_{y_i}}^{eq^n} \right). 
    \end{equation}
    Therefore, from \eqref{Ch_LB_Hineq_proof_a} and \eqref{Ch_LB_Hineq_proof_b}, we obtain
    \begin{eqnarray}
        H_q \left( f_{q_{x_i}}^{n+1} \right) - H_q \left(f_{q_{y_i}}^{n} \right) &\leq&  \left. \frac{\partial H_q}{\partial f_q} \right\rvert_{f_{q_{x_i}}^{n+1}} \left( f_{q_{x_i}}^{n+1}- f_{q_{y_i}}^{n} \right) \\
        &=& -\omega \left. \frac{\partial H_q}{\partial f_q} \right\rvert_{f_{q_{x_i}}^{n+1}} \left( f_{q_{x_i}}^{n+1} - f_{q_{y_i}}^{eq^n} \right) \\
        &\leq& - \omega \left( H_q \left( f_{q_{x_i}}^{n+1} \right) - H_q \left(f_{q_{y_i}}^{eq^n} \right) \right), \text{ since } \omega=\frac{\Delta t}{\epsilon}>0.  
    \end{eqnarray}
\end{proof}
\begin{remark}
\label{Ch_LB_Hineq_rem_semiimp}
    The following can be inferred from the above theorem:
    \begin{gather}
        H_q \left( f_{q_{x_i}}^{n+1} \right) \leq  \frac{1}{1+\omega} H_q \left(f_{q_{y_i}}^{n} \right) +  \frac{\omega}{1+\omega} H_q \left(f_{q_{y_i}}^{eq^n} \right), \\
        H_q \left( f_{q_{x_i}}^{n+1} \right) \leq  \left( 1-\Tilde{\omega} \right) H_q \left(f_{q_{y_i}}^{n} \right) +  \Tilde{\omega} H_q \left(f_{q_{y_i}}^{eq^n} \right). 
    \end{gather}
    Since $\sum_{q=1}^Q H_q\left(f_{q_{y_i}}^{eq^n}\right) \leq \sum_{q=1}^Q H_q\left(f_{q_{y_i}}^n\right)$ according to the definition of $H_q$, we obtain
    \begin{equation}
        \sum_{q=1}^Q H_q \left( f_{q_{x_i}}^{n+1} \right) \leq \sum_{q=1}^Q H_q \left(f_{q_{y_i}}^{n} \right). 
    \end{equation}
\end{remark}
Thus, for the LBE obtained from semi-implicit discretisation of the VKE, the H-inequality holds without enforcing any constraint on $\omega=\frac{\Delta t}{\epsilon}$. \\
The following remark \ref{Ch_LB_Hineq_rem_gen} presents H-inequality for general LBE, and the associated conditions. This has been presented particularly for explicit case in \cite{VK_Courtes}. 
\begin{remark}
\label{Ch_LB_Hineq_rem_gen}
    Consider the general LBE,
    \begin{eqnarray}
        f_{q_{x_i}}^{n+1} = \left( 1-\hat{\omega} \right) f_{q_{y_i}}^{n}  +  \hat{\omega} f_{q_{y_i}}^{eq^n}
    \end{eqnarray}
    with $\hat{\omega}=\omega$ (for explicit discretisation of VKE) and $\hat{\omega}=\Tilde{\omega}$ (for semi-implicit discretisation of VKE). Applying $H_q$ on this LBE, we obtain
    \begin{eqnarray}
        H_q \left( f_{q_{x_i}}^{n+1} \right) &=& H_q \left( \left( 1-\hat{\omega} \right) f_{q_{y_i}}^{n}  +  \hat{\omega} f_{q_{y_i}}^{eq^n} \right) \\
        &\leq& \left( 1-\hat{\omega} \right) H_q \left(f_{q_{y_i}}^{n}\right) + \hat{\omega} H_q\left(f_{q_{y_i}}^{eq^n} \right), \text{for } 0<\hat{\omega}\leq1 
    \end{eqnarray}
    Since $\sum_{q=1}^Q H_q\left(f_{q_{y_i}}^{eq^n}\right) \leq \sum_{q=1}^Q H_q\left(f_{q_{y_i}}^n\right)$, we obtain
    \begin{equation}
        \sum_{q=1}^Q H_q \left( f_{q_{x_i}}^{n+1} \right) \leq \sum_{q=1}^Q H_q \left(f_{q_{y_i}}^{n} \right). 
    \end{equation}
\end{remark}
Thus, the H-inequality holds for the general LBE if the constraint $0<\hat{\omega}\leq1$ is satisfied. It is to be noted that the H-inequality yields a stronger constraint on $\hat{\omega}$ than the positivity of numerical diffusion coefficient. \\
From the above remarks, the following can be inferred:  
\begin{itemize}
    \item For LBE obtained by explicit discretisation of VKE, $\hat{\omega}=\omega$. Hence, H-inequality holds corresponding to this LBE if $0<\omega=\frac{\Delta t}{\epsilon}\leq1$. It is to be noted that this constraint on $\omega$ is more restrictive than the constraint $0<\omega<2$ that enforces positivity of numerical diffusion coefficient. 
    \item For LBE obtained by semi-implicit discretisation of VKE, $\hat{\omega}=\Tilde{\omega}=\frac{\omega}{1+\omega}$. According to remark \ref{Ch_LB_Hineq_rem_gen}, H-inequality holds corresponding to this LBE if $0<\Tilde{\omega}=\frac{\omega}{1+\omega}\leq1$, and this is satisfied for all $\omega=\frac{\Delta t}{\epsilon}>0$. This also agrees with remark \ref{Ch_LB_Hineq_rem_semiimp} which states that H-inequality holds corresponding to this LBE for all $\omega=\frac{\Delta t}{\epsilon}>0$. Thus, the semi-implicit case of LBE is entropy-satisfying by construction.  
\end{itemize}

\subsection{Macroscopic finite difference form}
\label{Sec: LB_Macfd}
In this section, we show the macroscopic finite difference form of LBEs in \eqref{Ch_LB_LBEexplicit} (with $\hat{\omega}=\omega$) and \eqref{Ch_LB_LBEsemiimplicit} (with $\hat{\omega}=\Tilde{\omega}$). \\
We briefly provide some technicalities for clarity. LBE is evolved on a fixed uniform grid with spacing $\Delta x_d$ along direction $d$. At every time step, $\lambda_{d,n}$ is evaluated such that the sub-characteristic condition obtained by enforcing the positivity of numerical diffusion coefficient is satisfied. Thus, the discrete velocities can change with time step, and they are given by: $v_{q,n}^{d}=\left. m_q\right\rvert_d \lambda_{d,n}$, where $\left. m_q\right\rvert_d \in \mathbb{Z}$ is constant for direction $d$ and $q^{th}$ discrete velocity. The current time step is found by using $t_{n+1}-t_n:=\Delta t_n=\frac{\Delta x_d}{\lambda_{d,n}}$. Note that in addition to satisfying the sub-characteristic condition, $\lambda_{d,n}>\Delta x_d$ is essential for upper-bounding $\Delta t_n$ as $\Delta t_n<1$. Further, $\omega$ is kept constant for all time steps, and hence $\epsilon$ is allowed to depend on $n$ as $\Delta t_n$ depends on $n$.\\
For convenience, we consider one dimension ($D=1$) in the presentation of macroscopic finite difference form. Hence, the subscript and superscript $d$ indicating $d^{th}$ dimension can be ignored in all the variables. We consider the general LBE,
\begin{equation}
    f_q\left( x_i,t_n+\Delta t_n \right) = (1-\hat{\omega}) f_q \left( x_i-v_{q,n} \Delta t_n,t_n\right) +\hat{\omega} f^{eq}_q \left( U \left( x_i-v_{q,n}\Delta t_n,t_n\right) \right)
\end{equation}
with $\hat{\omega}=\omega$ and $\hat{\omega}=\Tilde{\omega}$ respectively for explicit and semi-implicit cases. For brevity, we introduce the following notations:\\ $f_{q_i}^{n+1}:=f_q\left( x_i,t_n+\Delta t_n \right)$, $f_{q_{i-m_q}}^n:=f_q \left( x_i-v_{q,n}\Delta t_n,t_n\right)$ and $f_{q_{i-m_q}}^{eq^n}:=f_q^{eq} \left( U \left( x_i-v_{q,n}\Delta t_n,t_n\right)\right)$. \\ We also utilise the splitting of $f_{q_i}^n$ as equilibrium and non-equilibrium parts: $f_{q_i}^{n}=f_{q_i}^{eq^n}+f_{q_i}^{neq^n}, \forall i,n$. Note here that we have absorbed $\epsilon$ of $\epsilon f_{q_i}^{neq^n}$ (refer \eqref{Ch_LB_DVBE_neq}) into $f_{q_i}^{neq^n}$ (\textit{i.e.,} $f_{q_i}^{neq^n}=\mathcal{O}(\epsilon)$) for convenience in presentation. Further, we also assume that $f_{q_i}^{0}=f_{q_i}^{eq^0}$ at the initial time. Thus, $f_{q_i}^{neq^0}=0$.  
\begin{theorem}
    The general LBE
    \begin{equation}
    \label{Ch_LB_Macfd_thm_genLBE}
        f_{q_i}^{n+1}=(1-\hat{\omega}) f_{q_{i-m_q}}^n + \hat{\omega} f_{q_{i-m_q}}^{eq^n}
    \end{equation}
    is equivalent to 
    \begin{equation}
    \label{Ch_LB_Macfd_thm_multstepLBE}
        f_{q_i}^{n+1}=\hat{\omega} \left( \sum_{k=0}^{N-1} (1-\hat{\omega})^k f_{q_{i-(k+1)m_q}}^{eq^{n-k}}  \right) + (1-\hat{\omega})^N f_{q_{i-(N+1)m_q}}^{eq^{n-N}}
    \end{equation}
    if $f_{q_{i-(N+1)m_q}}^{neq^{n-N}} = 0$. Here $N \in \mathbb{N}$.  
\end{theorem}
\begin{proof}
    Using $f_{q_{i-m_q}}^n=f_{q_{i-m_q}}^{eq^n}+f_{q_{i-m_q}}^{neq^n}$ in the general LBE \eqref{Ch_LB_Macfd_thm_genLBE}, we obtain
    \begin{equation}
    \label{Ch_LB_Macfd_proof_genLBE1}
        f_{q_i}^{n+1}= f_{q_{i-m_q}}^{eq^n} + (1-\hat{\omega}) f_{q_{i-m_q}}^{neq^n}.
    \end{equation}
    Using $f_{q_i}^{n+1}=f_{q_i}^{eq^{n+1}}+f_{q_i}^{neq^{n+1}}$ in the above equation yields 
    \begin{equation}
    \label{Ch_LB_Macfd_proof_genLBE2}
        f_{q_i}^{neq^{n+1}}=-f_{q_i}^{eq^{n+1}} + f_{q_{i-m_q}}^{eq^n} + (1-\hat{\omega}) f_{q_{i-m_q}}^{neq^n}.
    \end{equation}
    Inserting $f_{q_i}^{neq^{n+1}}$ from the above equation into $f_{q_{i-m_q}}^{neq^n}$ in \eqref{Ch_LB_Macfd_proof_genLBE1} by employing the transformation $n:=n'+1$, $i-m_q:=i'$, we get
    \begin{eqnarray}
        f_{q_i}^{n+1} &=& f_{q_{i-m_q}}^{eq^n} + (1-\hat{\omega}) \left( -f_{q_{i-m_q}}^{eq^{n}} + f_{q_{i-2m_q}}^{eq^{n-1}} + (1-\hat{\omega}) f_{q_{i-2m_q}}^{neq^{n-1}} \right), \\
        &=& \hat{\omega} f_{q_{i-m_q}}^{eq^{n}} + (1-\hat{\omega}) \left( f_{q_{i-2m_q}}^{eq^{n-1}} + (1-\hat{\omega}) f_{q_{i-2m_q}}^{neq^{n-1}} \right). 
    \end{eqnarray}
    Recursively inserting $f_{q_i}^{neq^{n+1}}$ from \eqref{Ch_LB_Macfd_proof_genLBE2} into the non-equilibrium term of above equation with the transformation $n-j:=n'+1$ and $i-(j+1)m_q:=i'$ where $j \in \{1,2,\dots, N-1 \}$, we get
    \begin{multline}
        f_{q_i}^{n+1}=\hat{\omega} \left( (1-\hat{\omega})^0 f_{q_{i-m_q}}^{eq^n} + (1-\hat{\omega})^1 f_{q_{i-2m_q}}^{eq^{n-1}} + \dots + (1-\hat{\omega})^{N-1} f_{q_{i-Nm_q}}^{eq^{n-(N-1)}} \right) \\ + (1-\hat{\omega})^{N} \left( f_{q_{i-(N+1)m_q}}^{eq^{n-N}} + (1-\hat{\omega}) f_{q_{i-(N+1)m_q}}^{neq^{n-N}} \right). 
    \end{multline}
    If $f_{q_{i-(N+1)m_q}}^{neq^{n-N}}=0$, then we obtain \eqref{Ch_LB_Macfd_thm_multstepLBE}. 
\end{proof} 

The above theorem depicts the multi-step nature of LBE by considering $t_{n-N}$ as the initial time. That is, $f_{q_i}^{n+1}$ depends on the values of the equilibrium function in neighboring grid points at all previous time steps starting from the initial time $t_{n-N}$. Note that $f_{q_{i-(N+1)m_q}}^{neq^{n-N}}=0$ as $f_{q_{i-(N+1)m_q}}^{n-N}=f_{q_{i-(N+1)m_q}}^{eq^{n-N}}$ is considered at the initial time.  \\
Summing \eqref{Ch_LB_Macfd_thm_multstepLBE} over $q$ with some form of equilibrium function discussed in section \ref{Ch_LB_Sec_Eqbm fun}, we obtain the macroscopic finite difference form. In this work, we consider the upwind D$\overline{d}$Q$(2\overline{d}+1)$ model (\textit{i.e.,} D$1$Q$3$ for one dimension). The equilibrium function for upwind D$1$Q$3$ model is,
\begin{gather}
    f^{eq^n}_{1_i}=\frac{G_i^{+^n}}{\lambda_n}\\
    f^{eq^n}_{2_i}=U_i^n-\frac{G_i^{+^n}+G_i^{-^n}}{\lambda_n} \\
    f^{eq^n}_{3_i}=\frac{G_i^{-^n}}{\lambda_n}
\end{gather}
and the corresponding velocities are $v_{1,n}=\lambda_n$, $v_{2,n}=0$ and $v_{3,n}=-\lambda_n$. Thus, $m_1=1$, $m_2=0$ and $m_3=-1$. 
\begin{remark}
    For $k \in \{0,1,\dots,N \}$, we have
    \begin{eqnarray}
        \sum_{q=1}^3 f_{q_{i-(k+1) m_q}}^{eq^{n-k}} &=& f_{1_{i-(k+1)}}^{eq^{n-k}}+f_{2_i}^{eq^{n-k}}+f_{3_{i+(k+1)}}^{eq^{n-k}} \\
        &=& \frac{G_{i-(k+1)}^{+^{n-k}}}{\lambda_{n-k}} + U_i^{n-k}-\frac{G_i^{+^{n-k}}+G_i^{-^{n-k}}}{\lambda_{n-k}} + \frac{G_{i+(k+1)}^{-^{n-k}}}{\lambda_{n-k}} \\
        &=& U_i^{n-k} - \frac{1}{\lambda_{n-k}} \left( \left( G_i^{+^{n-k}} -  G_{i-(k+1)}^{+^{n-k}} \right) - \left( G_{i+(k+1)}^{-^{n-k}} -  G_{i}^{-^{n-k}} \right) \right) \\
        &=& U_i^{n-k} - \frac{\Delta t_{n-k}}{\Delta x} \left( \left( G_i^{+^{n-k}} -  G_{i-(k+1)}^{+^{n-k}} \right) - \left( G_{i+(k+1)}^{-^{n-k}} -  G_{i}^{-^{n-k}} \right) \right).
    \end{eqnarray}
    Defining the notation 
    \begin{equation}
    \label{Ch_LB_Macfd_Uupdatenot}
        \mathcal{U}_{i,(k+1)}^{n-k+1}:=U_i^{n-k} - \frac{\Delta t_{n-k}}{\Delta x} \left( \left( G_i^{+^{n-k}} -  G_{i-(k+1)}^{+^{n-k}} \right) - \left( G_{i+(k+1)}^{-^{n-k}} -  G_{i}^{-^{n-k}} \right) \right),
    \end{equation} 
    $\sum_{q=1}^3$ \eqref{Ch_LB_Macfd_thm_multstepLBE} becomes
    \begin{equation}
    \label{Ch_LB_Macfdform_final}
        U_i^{n+1} = \hat{\omega} \left( \sum_{k=0}^{N-1} (1-\hat{\omega})^k \mathcal{U}_{i,(k+1)}^{n-k+1}  \right) + (1-\hat{\omega})^N \mathcal{U}_{i,(N+1)}^{n-N+1}.
    \end{equation}
\end{remark}
\noindent \eqref{Ch_LB_Macfdform_final} is the macroscopic finite difference form of the LBEs in \eqref{Ch_LB_LBEexplicit} (with $\hat{\omega}=\omega$) and \eqref{Ch_LB_LBEsemiimplicit} (with $\hat{\omega}=\Tilde{\omega}$).
\begin{remark}
    If $\hat{\omega}=1$, then $(1-\hat{\omega})^0 = 1$ and $(1-\hat{\omega})^k=0$ for $k \in \{1,2,\dots,N \}$. In this case, the macroscopic finite difference form \eqref{Ch_LB_Macfdform_final} becomes,
    \begin{equation}
        U_i^{n+1} = \mathcal{U}_{i,1}^{n+1} = U_i^{n} - \frac{\Delta t_{n}}{\Delta x} \left( \left( G_i^{+^{n}} -  G_{i-1}^{+^{n}} \right) - \left( G_{i+1}^{-^{n}} -  G_{i}^{-^{n}} \right) \right)
    \end{equation}
    which is an explicit (or forward) Euler upwind scheme for the hyperbolic system $\partial_t U + \partial_x G(U)=0$. 
\end{remark}
\textbf{Note:} For $k \in \{0,1,\dots,N \}$, $\mathcal{U}_{i,(k+1)}^{n-k+1}$ in \eqref{Ch_LB_Macfd_Uupdatenot} is an explicit (or forward) Euler upwind discretisation of the hyperbolic system $\partial_t U + \partial_x G(U)=0$, at time $t_{n-k+1}$ with grid spacing $(k+1)\Delta x$. Thus, \eqref{Ch_LB_Macfdform_final} which is the macroscopic finite difference form of LBE is simply a linear combination of upwind discretisations at varied time levels and grid spacings. 
\begin{remark}
\label{Ch_LB_remark_Numdiff}
    For $0<\hat{\omega}<1$, $(1-\hat{\omega})^k>0$ holds true for $k \in \{0,1,\dots,N \}$. Hence, in this case, numerical diffusion of the macroscopic finite difference form \eqref{Ch_LB_Macfdform_final} has positively weighted contributions from each $\mathcal{U}_{i,(k+1)}^{n-k+1}$. Thus, when $0<\hat{\omega}<1$, it is expected that the numerical diffusion increases with decrease in $\hat{\omega}$ while all the parameters remain frozen.  \\
    On the other hand, when $1<\hat{\omega}<2$, the sign of $(1-\hat{\omega})^k$ alternates with $k$. Therefore, numerical diffusion of the macroscopic finite difference form \eqref{Ch_LB_Macfdform_final} experiences alternately signed (with respect to $k$) weighted contributions from $\mathcal{U}_{i,(k+1)}^{n-k+1}$. 
\end{remark}
As a consequence, the minimum (over $\hat{\omega}$) numerical diffusion in LBE obtained by semi-implicit discretisation of VKE is larger than that in the explicit case.   
\subsection{Consistency}
In this section, we discuss the consistency of the macroscopic finite difference form \eqref{Ch_LB_Macfdform_final} with the hyperbolic system $\partial_t U + \partial_x G(U)=0$. 
\begin{theorem}
    Under suitable smoothness assumptions on all involved variables, the expression \eqref{Ch_LB_Macfd_Uupdatenot} becomes 
    \begin{eqnarray}
    \label{Ch_LB_Cons_thm1}
        \mathcal{U}_{i,(k+1)}^{n-k+1} = \left\{ \begin{matrix}
            U_i^n - \sum_{j=1}^k \Delta t_{n-j} \left. \partial_t U \right\rvert_i^n - (k+1) \Delta t_{n-k}  \partial_x \left. G \right\rvert_i^{n} & \text{for } k \in \{1,2,\dots,N \} \\
            U_i^n - (k+1) \Delta t_{n-k}  \partial_x \left. G \right\rvert_i^{n} & \text{for } k=0
        \end{matrix} \right.
    \end{eqnarray}
    upto $\mathcal{O}\left( k(k+1)\Delta x^2\right)$, if $\Delta t_m = \mathcal{O}(\Delta x) \ \forall m$.  
\end{theorem}
\begin{proof}
    Taylor expanding each term in $\mathcal{U}_{i,(k+1)}^{n-k+1}$:
    \begin{eqnarray}
        U_i^{n-k} &=& \left\{ \begin{matrix} U_i^n - \sum_{j=1}^k \Delta t_{n-j} \left. \partial_t U \right\rvert_i^n + \mathcal{O}\left( \left(k\Delta x\right)^2\right) & \text{for } k \in \{ 1,2,\dots,N \} \\
        U_i^n & \text{for } k=0
        \end{matrix}\right. \\
        && \text{since } \sum_{j=1}^k \Delta t_{n-j} = \mathcal{O}( k \Delta x) \left(\text{ as } \Delta t_m = \mathcal{O} (\Delta x),  \forall m\right) \nonumber
    \end{eqnarray}
    \begin{eqnarray*}
        \left( G_i^{+^{n-k}} -  G_{i-(k+1)}^{+^{n-k}} \right)  &=& (k+1) \Delta x \partial_x \left. G^+ \right\rvert_i^{n-k} + \mathcal{O}\left( \left((k+1)\Delta x\right)^2\right)\\
        \left( G_{i+(k+1)}^{-^{n-k}} -  G_{i}^{-^{n-k}} \right) &=& (k+1) \Delta x \partial_x \left. G^- \right\rvert_i^{n-k} + \mathcal{O}\left( \left((k+1)\Delta x\right)^2\right)
    \end{eqnarray*}
    \begin{eqnarray}
        \left( G_i^{+^{n-k}} -  G_{i-(k+1)}^{+^{n-k}} \right) - \left( G_{i+(k+1)}^{-^{n-k}} -  G_{i}^{-^{n-k}} \right) &=& (k+1) \Delta x \partial_x \left. \left( G^+ - G^- \right) \right\rvert_i^{n-k} \nonumber \\ &=& (k+1) \Delta x \partial_x \left. G \right\rvert_i^{n-k} \\
        && \text{upto } \mathcal{O}\left( \left((k+1)\Delta x\right)^2\right) \nonumber
    \end{eqnarray}
    \begin{equation}
        \partial_x \left. G \right\rvert_i^{n-k} = \left\{ \begin{matrix} 
        \partial_x \left. G \right\rvert_i^{n} - \sum_{j=1}^k \Delta t_{n-j} \partial_{tx}  \left. G \right\rvert_i^{n} + \mathcal{O}\left( \left(k\Delta x\right)^2\right) & \text{for } k \in \{1,2,\dots,N \} \\ \partial_x \left. G \right\rvert_i^{n} & \text{for } k =0 \end{matrix}  \right.
    \end{equation}
    For $k \in \{1,2,\dots,N \}$,
    \begin{eqnarray}
        (k+1)\Delta x \partial_x \left. G \right\rvert_i^{n-k} &\simeq & (k+1) \Delta x \partial_x \left. G \right\rvert_i^{n} - \sum_{j=1}^k \Delta t_{n-j} (k+1) \Delta x \partial_{tx}  \left. G \right\rvert_i^{n}  \nonumber \\
        &=& (k+1) \Delta x \partial_x \left. G \right\rvert_i^{n} + \mathcal{O}\left( k(k+1)\Delta x^2\right), 
    \end{eqnarray}
    since $\Delta t_m = \mathcal{O} (\Delta x), \forall m$. \\
    Thus, inserting the above expressions into \eqref{Ch_LB_Macfd_Uupdatenot}, we get \eqref{Ch_LB_Cons_thm1}. 
\end{proof}
\begin{remark}
    Taylor expanding $U_i^{n+1}$ about $U_i^n$, we get
    \begin{equation}
        U_i^{n+1} = U_i^n + \Delta t_n \left. \partial_t U \right\rvert_i^n + \mathcal{O}\left(\Delta t_n^2\right).
    \end{equation}
    Inserting \eqref{Ch_LB_Cons_thm1} and the above expression into \eqref{Ch_LB_Macfdform_final}, we obtain
    \begin{multline}
        U_i^n + \Delta t_n \left. \partial_t U \right\rvert_i^n = \hat{\omega} \biggl( (1-\hat{\omega})^0 \left( U_i^n -  \Delta t_{n} \partial_x \left. G \right\rvert_i^{n} \right) \biggr. \\ \left. +  \sum_{k=1}^{N-1} (1-\hat{\omega})^k \left( U_i^n - \sum_{j=1}^k \Delta t_{n-j} \left. \partial_t U \right\rvert_i^n - (k+1) \Delta t_{n-k}  \partial_x \left. G \right\rvert_i^{n} \right) \right) \\ + (1-\hat{\omega})^N \left( U_i^n - \sum_{j=1}^N \Delta t_{n-j} \left. \partial_t U \right\rvert_i^n - (N+1) \Delta t_{n-N}  \partial_x \left. G \right\rvert_i^{n} \right)
    \end{multline}
    upto $\mathcal{O}\left( N(N+1)\Delta x^2\right)$. Upon simplifying the above expression, we obtain the following upto $\mathcal{O}\left( N(N+1)\Delta x^2\right)$:
    \begin{multline}
    \label{Ch_LB_Cons_rem1}
        \left(1-\hat{\omega}\sum_{k=0}^{N-1}(1-\hat{\omega})^k  - (1-\hat{\omega})^N \right) U_i^n \\ + \left( \Delta t_n + \hat{\omega} \sum_{k=1}^{N-1} (1-\hat{\omega})^k \sum_{j=1}^k \Delta t_{n-j} +  (1-\hat{\omega})^N \sum_{j=1}^N \Delta t_{n-j} \right) \left. \partial_t U \right\rvert_i^n  \\ + \left( \hat{\omega} \sum_{k=0}^{N-1} (1-\hat{\omega})^k (k+1) \Delta t_{n-k} + (1-\hat{\omega})^N (N+1) \Delta t_{n-N} \right) \partial_x \left. G \right\rvert_i^{n} = 0. 
    \end{multline}
\end{remark}
\begin{remark}
    The coefficients of $U_i^n$ and $\left. \partial_t U \right\rvert_i^n$  in \eqref{Ch_LB_Cons_rem1} can be simplified as shown below:
    \begin{eqnarray}
    \label{Ch_LB_Cons_rem2_1}
        1-\hat{\omega}\sum_{k=0}^{N-1}(1-\hat{\omega})^k  - (1-\hat{\omega})^N  &=& 1- \hat{\omega} \left( \frac{1-\left( 1-\hat{\omega}\right)^N}{1-\left(1-\hat{\omega} \right)} \right) - (1-\hat{\omega})^N, \text{ for }\hat{\omega}\neq 0 \nonumber\\ &=& 0  
    \end{eqnarray}
    Since $1=\hat{\omega}\sum_{k=0}^{N-1}(1-\hat{\omega})^k + (1-\hat{\omega})^N$, we have $\Delta t_n=\left(\hat{\omega}\sum_{k=0}^{N-1}(1-\hat{\omega})^k + (1-\hat{\omega})^N\right)\Delta t_n$. Therefore, 
    \begin{equation}
    \label{Ch_LB_Cons_rem2_2}
        \Delta t_n + \hat{\omega} \sum_{k=1}^{N-1} (1-\hat{\omega})^k \sum_{j=1}^k \Delta t_{n-j} +  (1-\hat{\omega})^N \sum_{j=1}^N \Delta t_{n-j} = \hat{\omega}\sum_{k=0}^{N-1} (1-\hat{\omega})^k \sum_{j=0}^k \Delta t_{n-j} + (1-\hat{\omega})^N \sum_{j=0}^N \Delta t_{n-j}.
    \end{equation}
\end{remark}
Inserting \eqref{Ch_LB_Cons_rem2_1} and \eqref{Ch_LB_Cons_rem2_2} into \eqref{Ch_LB_Cons_rem1}, we obtain
\begin{multline}
\label{Ch_LB_Cons_main}
    \left( \hat{\omega}\sum_{k=0}^{N-1} (1-\hat{\omega})^k \sum_{j=0}^k \Delta t_{n-j} + (1-\hat{\omega})^N \sum_{j=0}^N \Delta t_{n-j} \right) \left. \partial_t U \right\rvert_i^n \\ + \left( \hat{\omega} \sum_{k=0}^{N-1} (1-\hat{\omega})^k (k+1) \Delta t_{n-k} + (1-\hat{\omega})^N (N+1) \Delta t_{n-N} \right) \partial_x \left. G \right\rvert_i^{n} = 0.
\end{multline}
upto $\mathcal{O}\left( N(N+1)\Delta x^2\right)$. 
\begin{remark}
\label{Ch_LB_Consistency_deltatconst}
    If $\Delta t_m = \Delta t, \forall m$, then \eqref{Ch_LB_Cons_main} becomes
    \begin{gather}
        \left( \hat{\omega} \sum_{k=0}^{N-1} (1-\hat{\omega})^k (k+1) \Delta t + (1-\hat{\omega})^N (N+1) \Delta t \right) \left(   \left. \partial_t U \right\rvert_i^n  + \partial_x \left. G \right\rvert_i^{n} \right) = \mathcal{O}\left( N(N+1)\Delta x^2\right), \nonumber \\
        \implies \left. \partial_t U \right\rvert_i^n  + \partial_x \left. G \right\rvert_i^{n} = \mathcal{O}\left( N\Delta x\right). 
    \end{gather}
    Thus, in this case, the macroscopic finite difference form of LBE is consistent with the hyperbolic system. 
\end{remark}
\begin{remark}
\label{Ch_LB_Consistency_omega1}
    If $\hat{\omega}=1$, then $(1-\hat{\omega})^0 = 1$ and $(1-\hat{\omega})^k=0$ for $k \in \{1,2,\dots,N \}$. Thus \eqref{Ch_LB_Cons_main} becomes
    \begin{gather}
        \Delta t_n \left(   \left. \partial_t U \right\rvert_i^n  + \partial_x \left. G \right\rvert_i^{n} \right) = \mathcal{O}\left( N(N+1)\Delta x^2\right), \nonumber \\
        \implies \left. \partial_t U \right\rvert_i^n  + \partial_x \left. G \right\rvert_i^{n} = \mathcal{O}\left( N (N+1) \Delta x\right). 
    \end{gather}
    Therefore, the macroscopic finite difference form of LBE is consistent with the hyperbolic system for this case too.
\end{remark}
Although the lattice Boltzmann algorithm is consistent for the two special cases: (i) constant time step size and (ii) $\hat{\omega}=1$, it can be seen from \eqref{Ch_LB_Cons_main} that consistency cannot be attained in the general case as $\sum_{j=0}^k \Delta t_{n-j} \neq (k+1) \Delta t_{n-k}$ for $k \in \{1,2,\dots,N \}$. However, one can choose constant $\Delta t$ such that the sub-characteristic condition holds for all time steps. In this way, the algorithm will be consistent with the hyperbolic system for the choice of the time-step satisfying the sub-characteristic condition.     

\subsection{Total Variation Boundedness}
The total variation boundedness (TVB) property of a numerical method for hyperbolic system ensures that the spatial variation remains bounded for all time steps. In this section, we discuss the TVB property of our lattice Boltzmann method by using its macroscopic finite difference form \eqref{Ch_LB_Macfdform_final}. This expression contains $\mathcal{U}_{i,(k+1)}^{n-k+1}$ for $k \in \{0,1,\dots,N \}$. For discussion of TVB property, we consider $\mathcal{U}_{i,(k+1)}^{n-k+1}$ derived by utilising upwind D$1$Q$3$ equilibrium function as in section \ref{Sec: LB_Macfd}. 
\begin{definition}
    The total variation of any variable $\theta$ defined on a lattice structure indexed by $i$ is given by,
\begin{equation*}
\mathbf{TV}(\theta) = \sum_i \left| \theta_{i+1}-\theta_i\right|
\end{equation*}
\end{definition}
\begin{theorem}
    Let $U_i^{n+1}$ given by \eqref{Ch_LB_Macfdform_final} be the macroscopic finite difference form. Then, its total variation satisfies 
    \begin{equation}
    \label{Ch_LB_TVD_thm1}
        \mathbf{TV} \left( U^{n+1} \right) \leq \left( \left| \hat{\omega} \right|  \sum_{k=0}^{N-1} \left| 1-\hat{\omega} \right|^k + \left| 1-\hat{\omega}  \right|^N   \right)  \mathrm{C}
    \end{equation}
    if $\mathbf{TV} \left( \mathcal{U}_{(k+1)}^{n-k+1} \right) \leq \mathrm{C}$, for $k \in \{0,1,\dots,N \}$. 
\end{theorem}
\begin{proof}
    Consider $U_i^{n+1}$ given by \eqref{Ch_LB_Macfdform_final}. Then, $U_{i+1}^{n+1}-U_i^{n+1}$ becomes
    \begin{equation*}
        U_{i+1}^{n+1}-U_i^{n+1} = \hat{\omega} \left( \sum_{k=0}^{N-1} (1-\hat{\omega})^k \left( \mathcal{U}_{i+1,(k+1)}^{n-k+1} - \mathcal{U}_{i,(k+1)}^{n-k+1} \right)  \right) + (1-\hat{\omega})^N \left( \mathcal{U}_{i+1,(N+1)}^{n-N+1} - \mathcal{U}_{i,(N+1)}^{n-N+1} \right). 
    \end{equation*}
    Then,
    \begin{eqnarray*}
        \left| U_{i+1}^{n+1}-U_i^{n+1} \right| &\leq& \left| \hat{\omega} \right| \sum_{k=0}^{N-1} \left| 1-\hat{\omega} \right|^k  \left|  \mathcal{U}_{i+1,(k+1)}^{n-k+1} - \mathcal{U}_{i,(k+1)}^{n-k+1} \right| + \left| 1-\hat{\omega}  \right|^N \left| \mathcal{U}_{i+1,(N+1)}^{n-N+1} - \mathcal{U}_{i,(N+1)}^{n-N+1} \right|, \\
        \implies \mathbf{TV} \left( U^{n+1} \right) &\leq& \left| \hat{\omega} \right| \sum_{k=0}^{N-1} \left| 1-\hat{\omega} \right|^k \mathbf{TV} \left( \mathcal{U}_{(k+1)}^{n-k+1} \right) + \left| 1-\hat{\omega}  \right|^N \mathbf{TV} \left( \mathcal{U}_{(N+1)}^{n-N+1} \right).
    \end{eqnarray*}
    Using $\mathbf{TV} \left( \mathcal{U}_{(k+1)}^{n-k+1} \right) \leq \mathrm{C}$ for $k \in \{0,1,\dots,N \}$ in the above expression, we get \eqref{Ch_LB_TVD_thm1}. 
\end{proof}
\begin{remark}
    If $0 < \hat{\omega} \leq 1$, then \eqref{Ch_LB_TVD_thm1} becomes
    \begin{eqnarray}
        \mathbf{TV} \left( U^{n+1} \right) &\leq& \left( \hat{\omega} \sum_{k=0}^{N-1} \left( 1-\hat{\omega} \right)^k + \left( 1-\hat{\omega}  \right)^N   \right)  \mathrm{C} \\
        &=& C, \text{ since }  \hat{\omega} \sum_{k=0}^{N-1} \left( 1-\hat{\omega} \right)^k + \left( 1-\hat{\omega}  \right)^N =1. 
    \end{eqnarray}
\end{remark}
Therefore, if the underlying difference scheme $\left( \mathcal{U}_{i,(k+1)}^{n-k+1} \right)$ due to the choice of equilibrium function is TVB, then the lattice Boltzmann method induced by it is also TVB (\textit{i.e.,} $\mathbf{TV} \left( U^{n+1} \right) \leq \mathrm{C}$) if $0 < \hat{\omega} \leq 1$. \\
Since upwind methods are TVB, $\mathbf{TV} \left( \mathcal{U}_{(k+1)}^{n-k+1} \right) \leq \mathrm{C}$ is true for the choice of upwind equilibrium function. Hence, the corresponding lattice Boltzmann method is also TVB if $0 < \hat{\omega} \leq 1$.

\subsection{Positivity}
Some of the variables of hyperbolic systems are positive for all time ({\em e.g.}, density and internal energy in Euler's system of gas dynamics, water height in shallow water system). Numerical schemes for such hyperbolic systems are expected to ensure the positivity of these variables. In this section, we show the positivity property of our lattice Boltzmann method by using its macroscopic finite difference form \eqref{Ch_LB_Macfdform_final}. $\mathcal{U}_{i,(k+1)}^{n-k+1}$ in this expression is derived by utilising upwind D$1$Q$3$ equilibrium function as in section \ref{Sec: LB_Macfd}. 
\begin{theorem}
    Let $U_i^{n+1}$ given by \eqref{Ch_LB_Macfdform_final} be the macroscopic finite difference form. If $\mathcal{U}_{i,(k+1)}^{n-k+1}$ is positive for $k \in \{0,1,\dots,N \}$ and $0 < \hat{\omega} \leq 1$, then $U_i^{n+1}$ is positive. 
\end{theorem}
\begin{proof}
    This is trivially seen from \eqref{Ch_LB_Macfdform_final}. 
\end{proof}

Therefore, if the underlying difference scheme $\left( \mathcal{U}_{i,(k+1)}^{n-k+1} \right)$ due to the choice of equilibrium function is positive, then the lattice Boltzmann method induced by it is also positive if $0 < \hat{\omega} \leq 1$. \\ \\
Thus, we discussed some properties of our LBEs. To conclude, the stability-related properties like H-inequality, total variation boundedness, and positivity are realisable if the stronger condition $0 < \hat{\omega} \leq 1$ is satisfied (naturally satisfied in the semi-implicit case) while small numerical diffusion is realisable for $\hat{\omega}>1$ (explicit case can be used in the interval $1 < \hat{\omega} < 2$ while ensuring positivity of numerical diffusion coefficient), depicting the trade-off between stability and accuracy.
\begin{remark}
    It is expected that the properties of LBM can be understood from its macroscopic finite difference form in \eqref{Ch_LB_Macfdform_final} by utilising the properties of corresponding underlying difference scheme $\mathcal{U}_{i,(k+1)}^{n-k+1}$ that occurs due to the choice of equilibrium functions. For instance, discrete conservation (with periodic boundary conditions) of LBM is evident if $\mathcal{U}_{i,(k+1)}^{n-k+1}$ satisfies discrete conservation with periodic boundary conditions.  
\end{remark}
Thus, in this section, novel discussions concerning LBEs derived by semi-implicit and explicit discretisations of VKE, on properties such as H-inequality, macroscopic finite difference form, consistency, total variation boundedness and positivity have been presented.  
\section{Hyperbolic conservation laws with source terms}
\label{Ch_LB_SecSource}
In this section, we extend our lattice Boltzmann method to hyperbolic conservation laws with source terms. Consider 
\begin{equation}
\label{Ch_LB_Hyp cons law_source}
    \partial_t U + \partial_{x_d} G^d(U)=S(U),
\end{equation}
where $S(U)$ is the source term.
\subsection{Vector-kinetic equation}
To approximate \eqref{Ch_LB_Hyp cons law_source}, consider the vector-kinetic equation
\begin{equation}
\label{Ch_LB_DVBE_source}
    \partial_t f_q + \partial_{x_d} \left( v^d_q f_q \right)= - \frac{1}{\epsilon} \left( f_q - f^{eq}_q(U) \right)+r(f_q).
\end{equation}
Summing \eqref{Ch_LB_DVBE_source} over all $q$, we get
\begin{equation}
    \partial_t \sum_{q=1}^Q f_q + \partial_{x_d} \sum_{q=1}^Q \left( v^d_q f_q \right)= - \frac{1}{\epsilon} \sum_{q=1}^Q \left( f_q - f^{eq}_q(U) \right) + \sum_{q=1}^Q r(f_q).
\end{equation}
If $\sum_{q=1}^Q f_q = \sum_{q=1}^Q f^{eq}_q = U$ and $\sum_{q=1}^Q r (f_q) =S(U)$, the above equation becomes
\begin{equation}
\label{Ch_LB_DVBE_sum_source}
    \partial_t U + \partial_{x_d} \sum_{q=1}^Q \left( v^d_q f_q \right)= S(U).
\end{equation}
In the limit $\epsilon \to 0$, we infer from \eqref{Ch_LB_DVBE_source} that $f_q \to f^{eq}_q(U)$. If $\sum_{q=1}^Q v^d_q f^{eq}_q = G^d(U)$, then \eqref{Ch_LB_DVBE_sum_source} becomes \eqref{Ch_LB_Hyp cons law_source} in the limit $\epsilon \to 0$. \\
Hereafter, we denote $r_q :=r(f_q)$ for convenience.

\subsection{Lattice Boltzmann equation}
As in section \ref{Ch_LB_SecLB}, $f_q$ in the collision term can be treated both explicitly and implicitly leading to LBEs with $\hat{\omega}=\omega$ and $\hat{\omega}=\Tilde{\omega}=\frac{\omega}{1+\omega}$ respectively. The source term $r_q$ is discretised in Crank-Nicolson fashion. Thus, the LBE becomes
\begin{multline}
\label{Ch_LB_LBE_source}
    f_q\left( \mathbf{x},t+\Delta t\right) = (1-\hat{\omega}) f_q \left( \mathbf{x}-\mathbf{v_q}\Delta t,t\right) +\hat{\omega} f^{eq}_q \left( U \left( \mathbf{x}-\mathbf{v_q}\Delta t,t\right) \right) \\ + \frac{\Delta t}{2} \left( r_q \left( \mathbf{x},t+\Delta t\right) + r_q \left( \mathbf{x}-\mathbf{v_q}\Delta t,t\right)   \right).
\end{multline}
The collision-streaming algorithm
\begin{eqnarray*}
    \text{Collision:} & F_q^* = (1-\hat{\omega}) f_q \left( \mathbf{x}-\mathbf{v_q}\Delta t,t\right) +\hat{\omega} f^{eq}_q \left( U \left( \mathbf{x}-\mathbf{v_q}\Delta t,t\right) \right) + \frac{\Delta t}{2} r_q \left( \mathbf{x}-\mathbf{v_q}\Delta t,t\right)\\
    \text{Streaming:} & F_q\left( \mathbf{x},t+\Delta t \right) = F_q^*\left( \mathbf{x}-\mathbf{v_q}\Delta t,t\right)
\end{eqnarray*}
can be used to numerically implement the LBEs. After finding $F_q\left( \mathbf{x},t+\Delta t \right) = f_q\left( \mathbf{x},t+\Delta t\right) - \frac{\Delta t}{2} r_q \left( \mathbf{x},t+\Delta t\right) $, we find $U\left( \mathbf{x},t+\Delta t\right)$ by solving 
\begin{equation}
\label{Ch_LB_nonlinearsolve_source}
    \sum_q F_q\left( \mathbf{x},t+\Delta t \right) = U \left( \mathbf{x},t+\Delta t\right) - \frac{\Delta t}{2} S(U\left( \mathbf{x},t+\Delta t\right))
\end{equation}
using a non-linear iterative solver (e.g., Newton's root finding method).

\subsection{Chapman-Enskog expansion}
The Chapman-Enskog expansion can be obtained by first Taylor expanding the LBE \eqref{Ch_LB_LBE_source} as,
\begin{equation}
\label{Ch_LB_CE_TE2_source}
\left( \partial_t + \mathbf{v_q} \cdot \nabla \right) f_q =  - \frac{\hat{\omega}}{\Delta t} \left(f_q-f^{eq}_q\right) + \frac{\hat{\omega}}{2} \left( \partial_t +  \mathbf{v_q} \cdot \nabla \right) \left(f_q-f_q^{eq}\right) + r_q +  \mathcal{O}(\Delta t^2).
\end{equation}  
Consider the perturbation expansions:
\begin{equation}
\label{Ch_LB_PE_source}
f_q=f_q^{eq} + \epsilon f_q^{(1)} + \epsilon^2 f_q^{(2)} + \dots; \ r_q= \epsilon r_q^{(1)} + \epsilon^2 r_q^{(2)} + \dots
\end{equation}
Since $\sum_{q=1}^Q f_q = \sum_{q=1}^Q f_q^{eq} = U$, we have $\sum_{q=1}^Q f_q^{(i)} = 0, \forall i \in \mathbb{N}$. Multiple scale expansion of derivatives of $f_q$ gives $\partial_t f_q = \left( \epsilon \partial_t^{(1)} + \epsilon^2 \partial_t^{(2)} + ... \right) f_q \ \text{and} \ \mathbf{v_q} \cdot \nabla f_q = \epsilon \mathbf{v_q} \cdot \nabla^{(1)} f_q$. \\
Using perturbation expansion of $f_q,r_q$ and multiple scale expansion of derivatives of $f_q$ in \eqref{Ch_LB_CE_TE2_source} and separating out $O(\epsilon)$ and $O(\epsilon^2)$ terms,
 \begin{eqnarray}
 \label{Ch_LB_CE_O(xi)_source}
 O(\epsilon): & \left( \partial_t^{(1)} + \mathbf{v_q} \cdot \nabla^{(1)}\right) f_q^{eq} = - \frac{\hat{\omega}}{\Delta t} f_q^{(1)} +  r_q^{(1)} \\
 \label{Ch_LB_CE_O(xi^2)_source}
 O(\epsilon^2): & \ \partial_t^{(2)} f_q^{eq} + \left( 1- \frac{\hat{\omega}}{2}\right)\left( \partial_t^{(1)} + \mathbf{v_q} \cdot \nabla^{(1)}\right) f_q^{(1)} =  - \frac{\hat{\omega}}{\Delta t} f_q^{(2)}  + r_q^{(2)}
 \end{eqnarray}
Zeroth moment $\left(\sum_{q=1}^Q\right)$ of $O(\epsilon)$ terms in \eqref{Ch_LB_CE_O(xi)_source} and $O(\epsilon^2)$ terms in \eqref{Ch_LB_CE_O(xi^2)_source} respectively give
 \begin{gather}
 \label{Ch_LB_CE_0thO(xi)_source}
 \partial_t^{(1)} U + \partial_{x_d}^{(1)} G^d(U) = \sum_{q=1}^Q  r_q^{(1)},  \\
 \label{Ch_LB_CE_0thO(xi^2)_source}
 \partial_t^{(2)} U + \left( 1- \frac{\hat{\omega}}{2}\right) \partial_{x_d}^{(1)} \left(\sum_{q=1}^Q v_q^{(d)}  f_q^{(1)}\right)=\sum_{q=1}^Q  r_q^{(2)}.  
 \end{gather}
From the first moment $\left(\sum_{q=1}^Q v^d_q\right)$ of $O(\epsilon)$ terms in \eqref{Ch_LB_CE_O(xi)_source}, we get 
\begin{equation}
\label{Ch_LB_CE_1stO(xi)_source}
\sum_{q=1}^Q v_q^d f_q^{(1)} = - \frac{\Delta t}{\hat{\omega}} \left( \partial_U G^d \left( -\partial_U G^i \partial_{x_i}^{(1)} U + \sum_{q=1}^Q  r_q^{(1)} \right) - \sum_{q=1}^Q v^d_q r^{(1)}_q + \partial_{x_i}^{(1)} \left( \sum_{q=1}^Q v_q^d v_q^i f_q^{eq} \right)\right)
\end{equation}
Recombining the zeroth moment equations of $O(\epsilon)$ in \eqref{Ch_LB_CE_0thO(xi)_source} and $O(\epsilon^2)$ in \eqref{Ch_LB_CE_0thO(xi^2)_source} and reversing the multiple scale expansions, we get 
\begin{multline}
\label{Ch_LB_CE_mPDE2_source}
\partial_t U + \partial_{x_d} G^d(U) = S(U) + \Delta t \left( \frac{1}{\hat{\omega}}-\frac{1}{2}\right) \left( \underbrace{\partial_{x_d} \left(  \partial_{x_i} \left( \sum_{q=1}^Q v_q^d v_q^i f_q^{eq} \right) - \partial_U G^d \partial_U G^i \partial_{x_i} U \right)}_{\text{Numerical Diffusion}}  \right. \\ \left. + \underbrace{\partial_{x_d} \left( \partial_U G^d \sum_{q=1}^Q  r_q - \sum_{q=1}^Q v^d_q r_q  \right)}_{\text{Spurious Numerical Convection}} \right).
\end{multline}

\subsection{Spurious Numerical Convection and modelling $r_q$}
The spurious numerical convection in \eqref{Ch_LB_CE_mPDE2_source} due to the discretisation of source term must by avoided in order to have a reliable numerical method. Therefore, we require $r_q$ to satisfy
\begin{equation}
\label{Ch_LB_vrcondn_source}
    \sum_{q=1}^Q v^d_q r_q = \partial_U G^d \sum_{q=1}^Q  r_q = \partial_U G^d S(U). 
\end{equation}
Thus, an $r_q$ that satisfies \eqref{Ch_LB_vrcondn_source} along with $\sum_{q=1}^Q r_q = S(U)$ is required. Note that these requirements are similar to those imposed on $f_q^{eq}$, and hence expressions similar to those in section \ref{Ch_LB_Sec_Eqbm fun} could be obtained for different models:
\begin{eqnarray*}
    \mathbf{\textbf{Classical } \textbf{D}1\textbf{Q}2:} & r_q=\frac{1}{2} S(U) - (-1)^q \frac{1}{2\lambda} \partial_U G^1 S(U), \text{ for } q \in \{1,2\} \\
    \mathbf{\textbf{D}1\textbf{Q}3:} & r_q=\frac{1}{3} S(U) + \left( \delta_{q1}-\delta_{q3} \right) \frac{1}{2\lambda} \partial_U G^1 S(U), \text{ for } q \in \{1,2,3\} \\
    \mathbf{\textbf{Upwind } \textbf{D}\overline{d}\textbf{Q}(2\overline{d}+1):} & r_q= \left\{ \begin{matrix}
        \frac{\partial_U G^{q+} S(U)}{\lambda_q}, & \text{for } q \in \{1,2,\dots,\overline{d} \} \\ 
        S(U) - \sum_{d=1}^{\overline{d}} \left( \frac{\left(\partial_U G^{d+}+\partial_U G^{d-}\right)S(U)}{\lambda_d} \right), & \text{for } q=\overline{d}+1 \\
        \frac{\partial_U G^{\left(q-\left(\overline{d}+1\right)\right)-} S(U)}{\lambda_{q-\left( \overline{d}+1 \right)}}, & \text{for } q \in \{\overline{d}+2,\dots,2\overline{d}+1 \}
    \end{matrix} \right.
\end{eqnarray*}
Thus, after finding $U \left( \mathbf{x},t+\Delta t\right)$ by solving \eqref{Ch_LB_nonlinearsolve_source}, we can find $f^{eq}_q\left(U\left( \mathbf{x},t+\Delta t \right)\right)$ (as discussed in section \ref{Ch_LB_Sec_Eqbm fun}) and $r_q\left(U \left( \mathbf{x},t+\Delta t \right) \right)$ (as discussed above), before proceeding with the next time step. 
\begin{remark}
\label{Ch_LB_source_remark}
    We modelled $r_q$ such that the spurious numerical convection due to discretisation of source term is nullified. This prevents the occurrence of spurious wave speeds and incorrect locations of discontinuities, commonly encountered in literature. Thus, balancing of convection and source terms which is a crucial problem in the finite volume framework, can be easily handled in the lattice Boltzmann framework.  Our strategy thus enforces the desired property of {\bf well-balancing} and, at the same time, takes care of stiffness of the source terms to a significant extent. \\
    Note that $r_q$ is of $\mathcal{O}(\epsilon)$ in \eqref{Ch_LB_PE_source}. Hence, our method and underlying removal of numerical convection works well for $S(U)=\sum_{q=1}^Q r_q = \mathcal{O}(\epsilon)$. 
\end{remark}

\section{D2Q9 model of lattice Boltzmann method}
\label{Ch_LB_SecD2Q9}
The equilibrium function in \eqref{Ch_LB_upweqbmfun} causes the underlying difference scheme $\left( \mathcal{U}_{i,(k+1)}^{n-k+1} \right)$ to result in pure upwinding along the coordinate directions. In this section, in addition to discrete velocities moving along coordinate directions, we also introduce discrete velocities moving along diagonal-to-coordinate directions. This enables the splitting of positive and negative fluxes even along diagonal-to-coordinate directions, thereby resulting in better multi-dimensional behavior. We consider two dimensions and a uniform lattice with equal grid spacing, with $\Delta x_1 = \Delta x_2 :=\Delta x$, in our presentation. 
\subsection{Equilibrium function}
We consider $9$ discrete velocities: $\mathbf{v_1}=[\lambda,0]$, $\mathbf{v_2}=[0,\lambda]$, $\mathbf{v_3}=[\lambda,\lambda]$, $\mathbf{v_4}=[-\lambda,\lambda]$, $\mathbf{v_5}=[0,0]$, $\mathbf{v_6}=[-\lambda,0]$, $\mathbf{v_7}=[0,-\lambda]$, $\mathbf{v_8}=[-\lambda,-\lambda]$, $\mathbf{v_9}=[\lambda,-\lambda]$, and the corresponding equilibrium functions: 
\begin{gather} 
f_1^{eq}=\frac{G^{\alpha+}}{\lambda}, \  f_2^{eq}=\frac{G^{\beta+}}{\lambda},\
f_3^{eq}=\frac{G^{\gamma+}}{\lambda}, \ f_4^{eq}=\frac{G^{\zeta+}}{\lambda}, \nonumber \\
\label{Ch_LB_D2Q9_Eqbm fn}
f_5^{eq}=U- \frac{1}{\lambda} \left( \left(G^{\alpha+} + G^{\beta+} +G^{\gamma+} + G^{\zeta+}\right) + \left(G^{\alpha-} + G^{\beta-} + G^{\gamma-} + G^{\zeta-}\right)\right),\\
f_6^{eq}=\frac{G^{\alpha-}}{\lambda}, \ f_7^{eq}=\frac{G^{\beta-}}{\lambda}, \
f_8^{eq}=\frac{G^{\gamma-}}{\lambda}, \ f_9^{eq}=\frac{G^{\zeta-}}{\lambda}. \nonumber
\end{gather}
Here $G^{l+} - G^{l-}=G^l$ for $l \in \mathcal{Z}=\{\alpha,\beta,\gamma,\zeta \}$.    
These equilibrium functions satisfy $\sum_{q=1}^{Q=9} f^{eq}_q = U$. In order to ensure $\sum_{q=1}^{Q=9} v^d_q f^{eq}_q = G^d(U)$, we need to satisfy the following requirements:
\begin{eqnarray}
    \sum_{q=1}^{Q=9} v^1_q f^{eq}_q = G^{\alpha} + G^{\gamma} - G^{\zeta} = G^1(U), \\
    \sum_{q=1}^{Q=9} v^2_q f^{eq}_q = G^{\beta} + G^{\gamma} + G^{\zeta} = G^2(U).
\end{eqnarray}
Thus, we have
\begin{equation}
G^{\gamma}=\frac{G^2+G^1}{2}-\frac{G^{\beta}+G^{\alpha}}{2} \text{ and } G^{\zeta}=\frac{G^2-G^1}{2}-\frac{G^{\beta}-G^{\alpha}}{2}, \ \forall G^{\alpha},G^{\beta}.
\end{equation}
In this setting, the underlying difference scheme corresponding to the equilibrium function \eqref{Ch_LB_D2Q9_Eqbm fn} is:
\begin{multline}
    \label{Ch_LB_D2Q9_unddiffsch}
        \mathcal{U}_{i,j,(k+1)}^{n-k+1}:=U_{i,j}^{n-k} - \frac{\Delta t_{n-k}}{\Delta x} \left( \left( G_{i,j}^{\alpha+^{n-k}} -  G_{i-(k+1),j}^{\alpha+^{n-k}} \right) - \left( G_{i+(k+1),j}^{\alpha-^{n-k}} -  G_{i,j}^{\alpha-^{n-k}} \right) \right) \\
        - \frac{\Delta t_{n-k}}{\Delta x} \left( \left( G_{i,j}^{\beta+^{n-k}} -  G_{i,j-(k+1)}^{\beta+^{n-k}} \right) - \left( G_{i,j+(k+1)}^{\beta-^{n-k}} -  G_{i,j}^{\beta-^{n-k}} \right) \right) \\
        - \frac{\Delta t_{n-k}}{\Delta x} \left( \left( G_{i,j}^{\gamma+^{n-k}} -  G_{i-(k+1),j-(k+1)}^{\gamma+^{n-k}} \right) - \left( G_{i+(k+1),j+(k+1)}^{\gamma-^{n-k}} -  G_{i,j}^{\gamma-^{n-k}} \right) \right) \\
        - \frac{\Delta t_{n-k}}{\Delta x} \left( \left( G_{i,j}^{\zeta+^{n-k}} -  G_{i+(k+1),j-(k+1)}^{\zeta+^{n-k}} \right) - \left( G_{i-(k+1),j+(k+1)}^{\zeta-^{n-k}} -  G_{i,j}^{\zeta-^{n-k}} \right) \right).
\end{multline} 
Further, the Chapman-Enskog expansion \eqref{Ch_LB_CE_mPDE2} corresponding to the equilibrium function \eqref{Ch_LB_D2Q9_Eqbm fn} becomes,
\begin{multline}
\label{Ch_LB_D2Q9_CE}
\partial_t U + \partial_{x_d} G^d(U) = \Delta t \left( \frac{1}{\hat{\omega}}-\frac{1}{2}\right) \\ \biggl( \partial_{x_1} \left( \lambda \partial_U \left(  G^{\alpha+} + G^{\alpha-} + G^{\gamma+} + G^{\gamma-} + G^{\zeta+} + G^{\zeta-} \right) - \left(\partial_U G^1\right)^2 \right) \partial_{x_1} U  \biggr. \\ \biggl. + \partial_{x_1}  \left( \lambda \partial_U \left( G^{\gamma+} + G^{\gamma-} - G^{\zeta+} - G^{\zeta-}  \right) - \partial_U G^1 \partial_U G^2 \right) \partial_{x_2} U \biggr. \\ \biggl. + \partial_{x_2}  \left( \lambda \partial_U \left( G^{\gamma+} + G^{\gamma-} - G^{\zeta+} - G^{\zeta-}  \right) - \partial_U G^2 \partial_U G^1 \right) \partial_{x_1} U \biggr. \\ + \biggl.
\partial_{x_2} \left( \lambda \partial_U \left( G^{\beta+} + G^{\beta-} + G^{\gamma+} + G^{\gamma-} + G^{\zeta+} + G^{\zeta-} \right) - \left(\partial_U G^2\right)^2 \right) \partial_{x_2} U \biggr).
\end{multline}  
Thus, in addition to upwinding along coordinate directions, this model allows upwinding even along diagonal-to-coordinate directions.   
\subsection{Boundary conditions}
In this sub-section, we present the expressions for $f_q$ corresponding to those specific $q$ that are unknown at the boundaries. At boundary, the macroscopic variables $U, G^{\alpha}, G^{\beta}, G^{\gamma} \text{ and }G^{\zeta}$ are known. From these, the split fluxes $G^{\alpha^{\pm}}, G^{\beta^{\pm}}, G^{\gamma^{\pm}}$ and $G^{\zeta^{\pm}}$ can be found. Using these split fluxes, equilibrium functions can be evaluated at the boundary. Thus, by taking $f_q^{neq}=f_q-f_q^{eq} \ \forall q \in \{1,2,..,9 \}$, it can be inferred from the definition of conserved moment $\sum_{n=1}^9 f_q = \sum_{q=1}^9 f_q^{eq}=U$ that, $\sum_{q=1}^9 f_q^{neq}=0$. 

 \begin{figure}[!h]
 \centering
 \begin{subfigure}[b]{0.13\textwidth}
 \centering
\includegraphics[width=\textwidth]{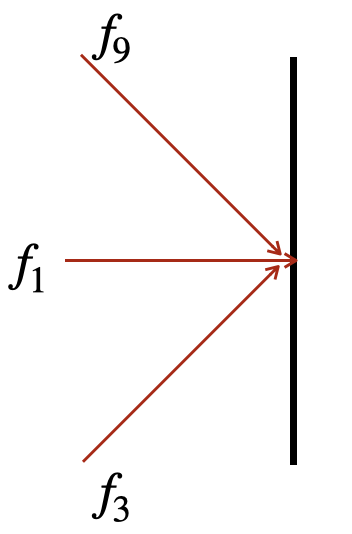}
\caption{Left}
\label{Ch_LB_Fig:D2Q9_left}
\end{subfigure}
\hfill
 \begin{subfigure}[b]{0.13\textwidth}
\centering
\includegraphics[width=\textwidth]{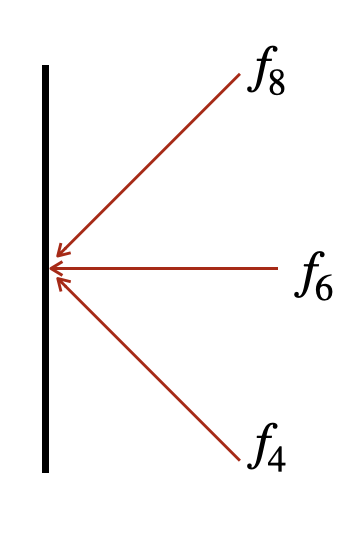}
\caption{Right}
\label{Ch_LB_Fig:D2Q9_right}
\end{subfigure}
\hfill
\begin{subfigure}[b]{0.2\textwidth}
\centering
\includegraphics[width=\textwidth]{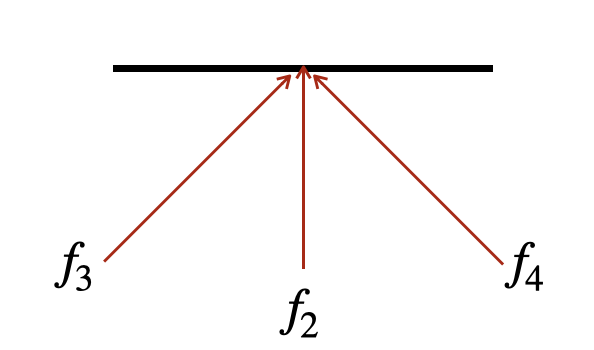}
\caption{Bottom}
\label{Ch_LB_Fig:D2Q9_bottom}
\end{subfigure}
\hfill
\begin{subfigure}[b]{0.2\textwidth}
\centering
\includegraphics[width=\textwidth]{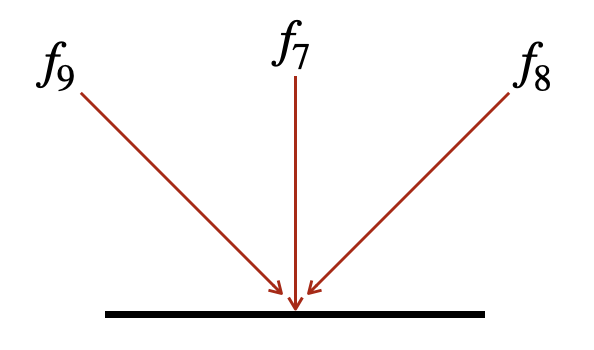}
\caption{Top}
\label{Ch_LB_Fig:D2Q9_top}
\end{subfigure}
\caption{\centering Boundary conditions (Black lines indicate boundaries; red arrows indicate unknown functions at each boundary)}
\end{figure}

\subsubsection{Left boundary}
At any point on left boundary, $f_2, f_4, f_5,f_6,f_7 \text{ and } f_8$ are known from the computational domain, as these functions from neighbouring points (in the computational domain) hop to points on left boundary. Let $\mathcal{I}$ be the set of these known functions. The unknowns at left boundary are $f_1, f_3 \text{ and } f_9$ (as shown in figure \ref{Ch_LB_Fig:D2Q9_left}), as these functions must come from the outside of computational domain to left boundary, and let $\mathcal{J}$ be the set of these unknown functions. Since $f_q^{eq}$ can be evaluated $\forall q \in \{1,2,..,9 \}$ and $f_q$ is known $\forall q \in \mathcal{I}$, $f_q^{neq}=f_q-f_q^{eq}$ can be found $\forall q \in \mathcal{I} $ (as $\mathcal{I} \subset \{1,2,..,9 \}$). Then $f_q^{neq}, \forall q \in \mathcal{J}$ can be written as, 
\begin{eqnarray}
f_3^{neq}=-f_8^{neq}-\frac{f_2^{neq}+f_5^{neq}+f_7^{neq}}{3} \\
f_1^{neq}=-f_6^{neq}-\frac{f_2^{neq}+f_5^{neq}+f_7^{neq}}{3} \\
f_9^{neq}=-f_4^{neq}-\frac{f_2^{neq}+f_5^{neq}+f_7^{neq}}{3}
\end{eqnarray} 
satisfying $\sum_{q=1}^9 f_q^{neq}=0$. Now, $f_q=f_q^{eq}+f_q^{neq} \  \forall q \in \mathcal{J}$ can be found to be,
\begin{eqnarray}
f_3=\frac{G^{\gamma+}+G^{\gamma-}}{\lambda}+\frac{U}{3}-\frac{1}{3\lambda} \sum\limits_{z\in \mathcal{Z}, z \neq \beta} \left(G^{z+}+G^{z-}\right)-f_8-\frac{f_2+f_5+f_7}{3} \\
f_1=\frac{G^{\alpha+}+G^{\alpha-}}{\lambda}+\frac{U}{3}-\frac{1}{3\lambda} \sum\limits_{z\in \mathcal{Z}, z \neq \beta} \left(G^{z+}+G^{z-}\right)-f_6-\frac{f_2+f_5+f_7}{3} \\
f_9=\frac{G^{\zeta+}+G^{\zeta-}}{\lambda}+\frac{U}{3}-\frac{1}{3\lambda} \sum\limits_{z\in \mathcal{Z}, z \neq \beta} \left(G^{z+}+G^{z-}\right)-f_4-\frac{f_2+f_5+f_7}{3}
\end{eqnarray}

\subsubsection{Right boundary}
By following the same procedure of obtaining left boundary conditions, the unknown functions at right boundary (as shown in figure \ref{Ch_LB_Fig:D2Q9_right}) can be found as, 
\begin{eqnarray}
f_4=\frac{G^{\zeta+}+G^{\zeta-}}{\lambda}+\frac{U}{3}-\frac{1}{3\lambda} \sum\limits_{z\in \mathcal{Z}, z \neq \beta} \left(G^{z+}+G^{z-}\right)-f_9-\frac{f_2+f_5+f_7}{3} \\
f_6=\frac{G^{\alpha+}+G^{\alpha-}}{\lambda}+\frac{U}{3}-\frac{1}{3\lambda} \sum\limits_{z\in \mathcal{Z}, z \neq \beta} \left(G^{z+}+G^{z-}\right)-f_1-\frac{f_2+f_5+f_7}{3} \\
f_8=\frac{G^{\gamma+}+G^{\gamma-}}{\lambda}+\frac{U}{3}-\frac{1}{3\lambda} \sum\limits_{z\in \mathcal{Z}, z \neq \beta} \left(G^{z+}+G^{z-}\right)-f_3-\frac{f_2+f_5+f_7}{3}
\end{eqnarray}

\subsubsection{Bottom boundary}
The unknown functions at bottom boundary (as shown in figure \ref{Ch_LB_Fig:D2Q9_bottom}) can be found as, 
\begin{eqnarray}
f_3=\frac{G^{\gamma+}+G^{\gamma-}}{\lambda}+\frac{U}{3}-\frac{1}{3\lambda} \sum\limits_{z\in \mathcal{Z}, z \neq \alpha} \left(G^{z+}+G^{z-}\right)-f_8-\frac{f_1+f_5+f_6}{3} \\
f_2=\frac{G^{\beta+}+G^{\beta-}}{\lambda}+\frac{U}{3}-\frac{1}{3\lambda} \sum\limits_{z\in \mathcal{Z}, z \neq \alpha} \left(G^{z+}+G^{z-}\right)-f_7-\frac{f_1+f_5+f_6}{3} \\
f_4=\frac{G^{\zeta+}+G^{\zeta-}}{\lambda}+\frac{U}{3}-\frac{1}{3\lambda} \sum\limits_{z\in \mathcal{Z}, z \neq \alpha} \left(G^{z+}+G^{z-}\right)-f_9-\frac{f_1+f_5+f_6}{3}
\end{eqnarray}

\subsubsection{Top boundary}
The unknown functions at top boundary (as shown in figure \ref{Ch_LB_Fig:D2Q9_top}) can be found as, 
\begin{eqnarray}
f_9=\frac{G^{\zeta+}+G^{\zeta-}}{\lambda}+\frac{U}{3}-\frac{1}{3\lambda} \sum\limits_{z\in \mathcal{Z}, z \neq \alpha} \left(G^{z+}+G^{z-}\right)-f_4-\frac{f_1+f_5+f_6}{3} \\
f_7=\frac{G^{\beta+}+G^{\beta-}}{\lambda}+\frac{U}{3}-\frac{1}{3\lambda} \sum\limits_{z\in \mathcal{Z}, z \neq \alpha} \left(G^{z+}+G^{z-}\right)-f_2-\frac{f_1+f_5+f_6}{3} \\
f_8=\frac{G^{\gamma+}+G^{\gamma-}}{\lambda}+\frac{U}{3}-\frac{1}{3\lambda} \sum\limits_{z\in \mathcal{Z}, z \neq \alpha} \left(G^{z+}+G^{z-}\right)-f_3-\frac{f_1+f_5+f_6}{3}
\end{eqnarray}

\subsubsection{Bottom-left corner}
At bottom left corner, the known equilibrium functions are $f_7, f_8, f_5 \text{ and } f_6$. The unknown equilibrium functions are $f_1, f_3, f_2, f_4 \text{ and }f_9$. Since $f_4 \text{ and } f_9$ do not enter or leave the computational domain, evaluation of them is not needed. Hence, it can be assumed that $f_9^{neq}+f_4^{neq}+f_5^{neq}=0$. Then $f_q^{neq}$ for other unknown equilibrium distribution functions can be written as, 
\begin{eqnarray}
f_1^{neq}=-f_6^{neq} \\
f_3^{neq}=-f_8^{neq} \\
f_2^{neq}=-f_7^{neq}
\end{eqnarray}
satisfying $\sum_{q=1}^9 f_q^{neq}=0$. Now, $f_q=f_q^{eq}+f_q^{neq}$ can be found to be, 
\begin{eqnarray}
f_1=\frac{G^{\alpha+}+G^{\alpha-}}{\lambda}-f_6 \\
f_3=\frac{G^{\gamma+}+G^{\gamma-}}{\lambda}-f_8 \\
f_2=\frac{G^{\beta+}+G^{\beta-}}{\lambda}-f_7
\end{eqnarray}

 \begin{figure}[!h]
 \centering
 \begin{subfigure}[b]{0.2\textwidth}
 \centering
\includegraphics[width=\textwidth]{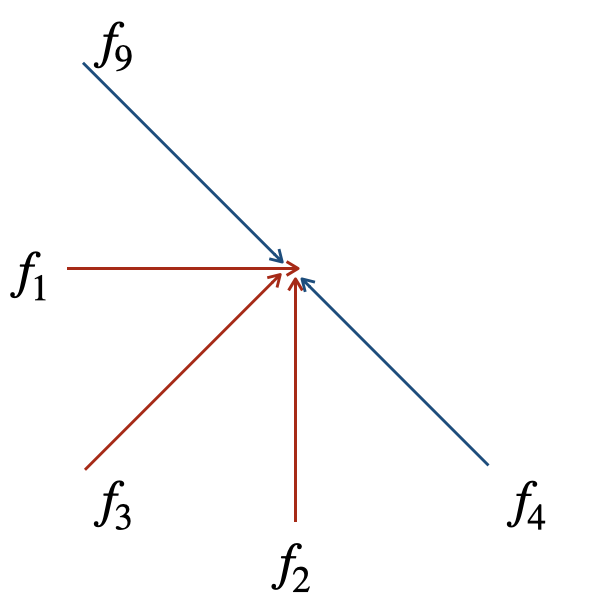}
\caption{Bottom-left}
\label{Ch_LB_Fig:D2Q9_bottomleft}
\end{subfigure}
\hfill
 \begin{subfigure}[b]{0.2\textwidth}
\centering
\includegraphics[width=\textwidth]{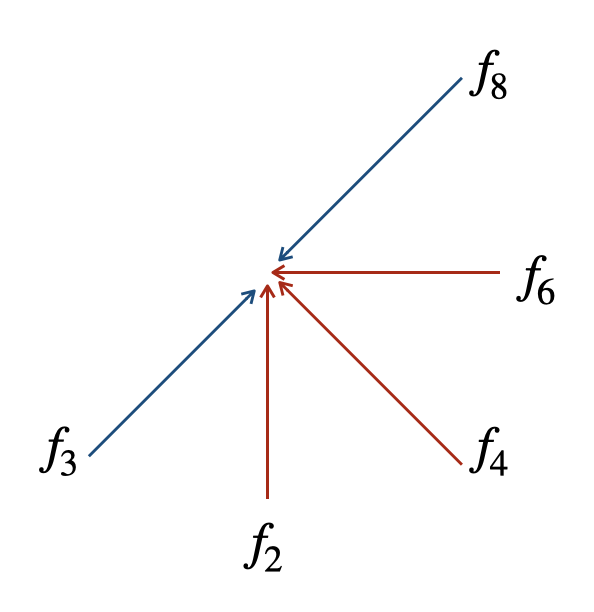}
\caption{Bottom-right}
\label{Ch_LB_Fig:D2Q9_bottomright}
\end{subfigure}
\hfill
\begin{subfigure}[b]{0.2\textwidth}
\centering
\includegraphics[width=\textwidth]{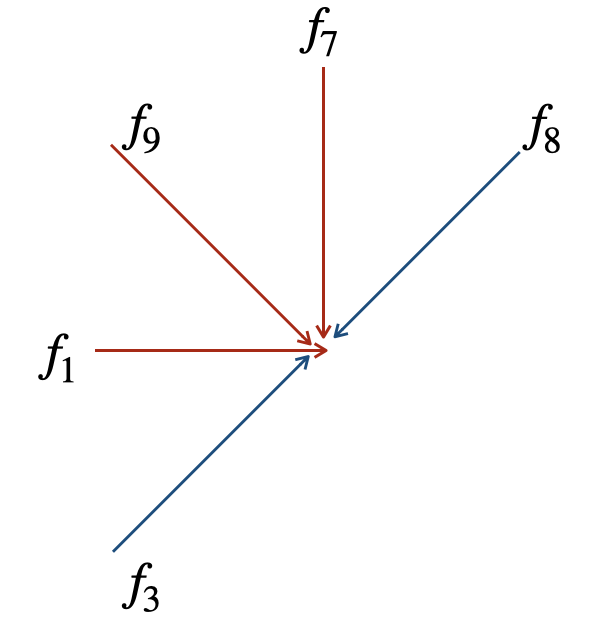}
\caption{Top-left}
\label{Ch_LB_Fig:D2Q9_topleft}
\end{subfigure}
\hfill
\begin{subfigure}[b]{0.2\textwidth}
\centering
\includegraphics[width=\textwidth]{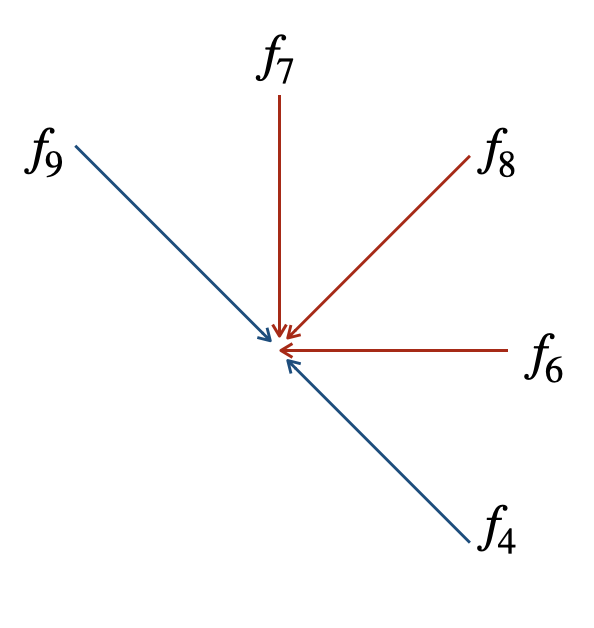}
\caption{Top-right}
\label{Ch_LB_Fig:D2Q9_topright}
\end{subfigure}
\caption{\centering Corner conditions (Red arrows indicate unknown functions that are evaluated; Blue arrows indicate unknown functions that are not evaluated)}
\end{figure}

\subsubsection{Bottom-right corner}
By following the same procedure for obtaining bottom-left corner conditions, the bottom-right corner conditions (as shown in figure \ref{Ch_LB_Fig:D2Q9_bottomright}) are found to be,
\begin{eqnarray}
f_2=\frac{G^{\beta+}+G^{\beta-}}{\lambda}-f_7 \\
f_4=\frac{G^{\zeta+}+G^{\zeta-}}{\lambda}-f_9 \\
f_6=\frac{G^{\alpha+}+G^{\alpha-}}{\lambda}-f_1
\end{eqnarray}

\subsubsection{Top-left corner}
The top-left corner conditions (as shown in figure \ref{Ch_LB_Fig:D2Q9_topleft}) are,
\begin{eqnarray}
f_1=\frac{G^{\alpha+}+G^{\alpha-}}{\lambda}-f_6 \\
f_9=\frac{G^{\zeta+}+G^{\zeta-}}{\lambda}-f_4 \\
f_7=\frac{G^{\beta+}+G^{\beta-}}{\lambda}-f_2
\end{eqnarray}

\subsubsection{Top-right corner}
The top-right corner conditions (as shown in figure \ref{Ch_LB_Fig:D2Q9_topright}) are,
\begin{eqnarray}
f_6=\frac{G^{\alpha+}+G^{\alpha-}}{\lambda}-f_1 \\
f_7=\frac{G^{\beta+}+G^{\beta-}}{\lambda}-f_2 \\
f_8=\frac{G^{\gamma+}+G^{\gamma-}}{\lambda}-f_3
\end{eqnarray}

\section{Numerical results}
\label{Ch_LB_SecNum}
In this section, we present the numerical validation of our lattice Boltzmann methods (LBM) discussed in the previous sections. Firstly, we depict the influence of $\omega$ on numerical diffusion and order of accuracy. Then, we numerically validate our LBM for hyperbolic conservation laws with source terms, and D$2$Q$9$ model of LBM. For all the cases, the numerical results are obtained by using LBE derived by explicit discretisation of VKE. Due to the algorithmic similarity of LBEs derived by explicit and semi-implicit discretisation of VKE, the numerical results obtained by semi-implicit case for $0<\Tilde{\omega}<1$ are same as that obtained by explicit case for $0<\omega<1$. Hence we only present the numerical validation of explicit case with larger interval $0<\omega<2$.  
\subsection{Sinusoidal initial condition}
The domain of the problem is $[0,1] \subset \mathbb{R}$. We consider inviscid Burgers' equation with flux function as $G^1(U)=\frac{1}{2}U^2$. The initial condition is $U(x_1,0)=\sin (2\pi x_1)$. An LBM with upwind D$1$Q$3$ equilibrium functions is utilised to obtain the numerical solution. $\lambda_1=\frac{\Delta x_1}{\Delta t}$ is chosen such that the sub-characteristic condition in \eqref{Ch_LB_Upwindeqbm_subcharcondn} (which simplifies in this case as $\lambda_1 \geq \sup\limits_{i \in \Omega_g} |U_i|$, where $\Omega_g$ is the set of grid points) is satisfied. Since we expect the numerical solution to be bounded between $-1$ and $1$ for all times, we choose $\lambda_1=1$, and fix $\Delta t = \frac{\Delta x_1}{\lambda_1}$ for all time steps in order to have a consistent discretisation of the inviscid Burgers' equation (as discussed in remark \ref{Ch_LB_Consistency_deltatconst}). Further, we consider different values for $\omega=\frac{\Delta t}{\epsilon}$ such as, $\omega=0.1,0.6,1.0,1.4,1.9$ and compare their numerical diffusion by freezing all the other parameters. We also consider discretisation of the domain with different number of grid points $N$ such as, $N=41,81,161,321$ in order to study the order of convergence. The reference solution utilised in finding the $L_2$ error norm is obtained by evaluating the method of characteristics solution with a tolerance of $10^{-15}$.  \\
Tables \ref{tab:Ch_LB_NR_Sinusoid_1} and \ref{tab:Ch_LB_NR_Sinusoid_2} show the $L_2$ error norms and convergence orders for different values of $\omega$ at time $T=\frac{0.1}{2\pi}$ while the solution is still smooth. It is seen from the tables that for each fixed value of $N$, $L_2$ error norm of the numerical solution increases with decrease in $\omega$, validating the remark \ref{Ch_LB_remark_Numdiff}. Further, although only first order of accuracy is expected according to Chapman-Enskog expansion \eqref{Ch_LB_Eqbmupwind_CEexp}, we observe more than second order accuracy for large values of $\omega$. This increase in order of accuracy for large values of $\omega$ can be attributed to the smaller numerical diffusion for $\omega > 1$ when compared to $\omega < 1$, as mentioned in remark \ref{Ch_LB_remark_Numdiff}. We also observe that $\mathcal{O}(L_2)$ corresponding to a fixed $N$ increases with increase in $\omega$.  

\begin{table}[tbhp]
\begin{center}
\begin{tabular}{|m{0.5cm}|m{1.45cm}|m{1.7cm}|m{1.9cm}|m{1.7cm}|m{1.9cm}|m{1.7cm}|m{1.9cm}|}
\hline
\centering $N$ & \centering $\Delta x_1$ & \centering $L_2, \ \omega=1.9 $ &  $\mathcal{O}(L_2), \ \omega=1.9$ & \centering $L_2, \ \omega=1.4$ & $\mathcal{O}(L_2), \  \omega=1.4$ & \centering $L_2, \ \omega=1.0$ &  $\mathcal{O}(L_2), \ \omega=1.0$ \\
\hline
\centering 41 & \centering 0.025 & \centering 0.000597 &  - & \centering 0.000597 &  - & \centering 0.000597 &  - \\
\centering 81 & \centering 0.0125 & \centering 9.68 $\times 10^{-5}$ &  2.626 & \centering 0.000158 &  1.915 & \centering 0.000230 &  1.380\\
\centering 161 & \centering 0.00625 & \centering 2.14 $\times 10^{-5}$ &  2.175 & \centering 3.88$\times 10^{-5}$ &  2.032  & \centering 6.41 $\times 10^{-5}$ &  1.841 \\
\centering 321 & \centering 0.003125 & \centering 3.20 $\times 10^{-6}$ &  2.744 & \centering 1.20 $\times 10^{-5}$ &  1.690  & \centering 2.33 $\times 10^{-5}$ &  1.460 \\
\hline
\end{tabular}
\caption{\centering Sinusoidal initial condition at $T=\frac{0.1}{2\pi}$ for $\omega=1.9,1.4,1.0$} 
\label{tab:Ch_LB_NR_Sinusoid_1}
\end{center}
\end{table}

\begin{table}[tbhp]
\begin{center}
\begin{tabular}{|m{0.5cm}|m{1.45cm}|m{1.7cm}|m{1.9cm}|m{1.7cm}|m{1.9cm}|}
\hline
\centering $N$ & \centering $\Delta x_1$ & \centering $L_2, \ \omega=0.6$ & $\mathcal{O}(L_2), \ \omega=0.6$ & \centering $L_2,\ \omega=0.1$ & $\mathcal{O}(L_2), \ \omega=0.1$ \\
\hline
\centering 41 & \centering 0.025 & \centering 0.000597 &  - & \centering 0.000597 &  -  \\
\centering 81 & \centering 0.0125 & \centering 0.000306 &  0.965 & \centering 0.000405 &  0.562 \\
\centering 161 & \centering 0.00625 & \centering 0.000100 &  1.611 & \centering 0.000161 &  1.325  \\
\centering 321 & \centering 0.003125 & \centering 4.38 $\times 10^{-5}$ &  1.194 & \centering 0.000103 &  0.644  \\
\hline
\end{tabular}
\caption{\centering Sinusoidal initial condition at $T=\frac{0.1}{2\pi}$ for $\omega=0.6,0.1$} 
\label{tab:Ch_LB_NR_Sinusoid_2}
\end{center}
\end{table}

\subsection{LBM for hyperbolic conservation laws with source terms}
The governing equation is of the form \eqref{Ch_LB_Hyp cons law_source} with p=1 (scalar conservation law). We show that our scheme captures the discontinuities at correct locations due to the nullification of spurious numerical convection by our choice of $r_q$. Further, since $r_q=\mathcal{O}(\epsilon)$ is essential for such a possibility of nullification as mentioned in remark \ref{Ch_LB_source_remark}, the numerical results with correct locations of discontinuities are presented whenever $S(U)=\mathcal{O}(\epsilon)$. 
\subsubsection{One dimensional discontinuity}
\label{Ch_LB_Numres_source_1D}
This is the test problem used by LeVeque and Yee \cite{LeVequeYee1990} to understand the cause for incorrectness in speeds of discontinuities for stiff source terms. The domain is $[0,1] \subset \mathbb{R}$, and is split up into $50$ evenly spaced grid points. For this problem, $G^1(U)=U$ and $S(U)=-\mu U (U-1) (U-\frac{1}{2})$. Initial conditions are: 
\begin{equation*}
    U(x_1,0)=\left\{ \begin{matrix}
        1 & \text{for } x_1 \leq 0.3 \\
        0 & \text{for } x_1 > 0.3
    \end{matrix} \right.. 
\end{equation*}
An LBM with upwind D$1$Q$3$ form for $f_q^{eq}$ and $r_q$ is utilised to obtain the numerical solution. $\lambda_1=\frac{\Delta x_1}{\Delta t}$ is chosen such that the sub-characteristic condition in \eqref{Ch_LB_Upwindeqbm_subcharcondn} (which simplifies in this case as $\lambda_1 \geq 1$) is satisfied. In particular, we use $\lambda_1=1$, and this incidentally results in numerical solution being the same as method of characteristics solution (even in smooth regions) since the wave-speed in the problem is also $1$. Therefore, in addition to capturing discontinuities at correct locations (due to our choice of $r_q$), the solution is also exact in smooth regions. Further, the time step is chosen as $\Delta t = \frac{\Delta x_1}{\lambda_1}$. We also consider $\omega=1$ for the simulation of this problem and this ensures consistency with the governing equation irrespective of the choice of $\Delta t$ (as discussed in remark \ref{Ch_LB_Consistency_omega1}). \\ 
\begin{figure}[!h]
\centering
\begin{subfigure}[b]{0.24\textwidth}
\centering
\includegraphics[width=\textwidth]{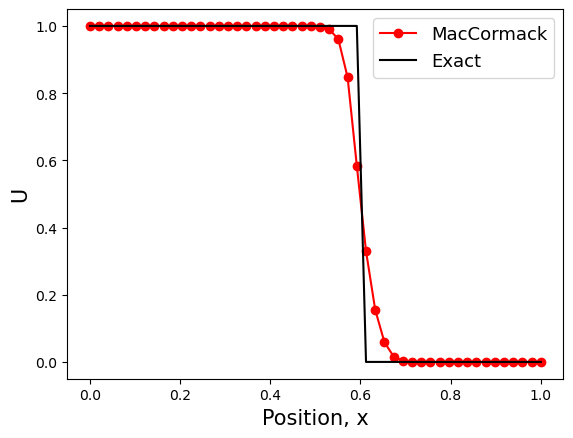}
\caption{$\mu=1$}
\label{Fig:Ch_LB_1DLY_MC_mu1}
\end{subfigure}
\hfill
\begin{subfigure}[b]{0.24\textwidth}
\centering
\includegraphics[width=\textwidth]{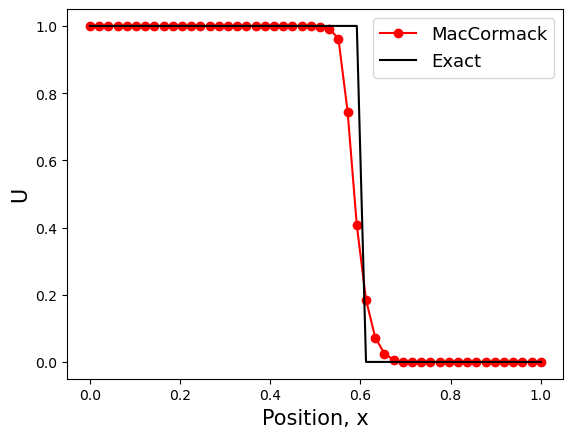}
\caption{$\mu=10$}
\label{Fig:Ch_LB_1DLY_MC_mu10}
\end{subfigure}
\hfill
\begin{subfigure}[b]{0.24\textwidth}
\centering
\includegraphics[width=\textwidth]{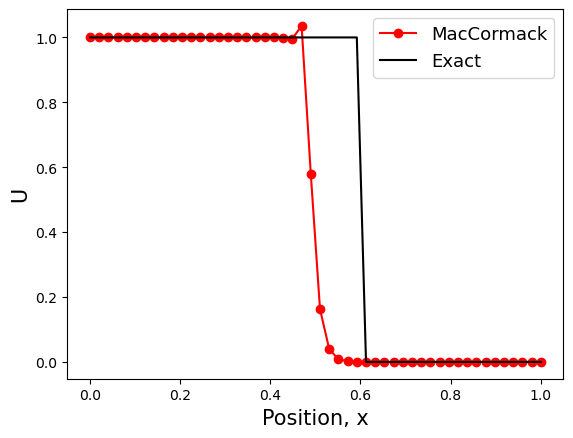}
\caption{$\mu=100$}
\label{Fig:Ch_LB_1DLY_MC_mu100}
\end{subfigure}
\hfill
\begin{subfigure}[b]{0.24\textwidth}
\centering
\includegraphics[width=\textwidth]{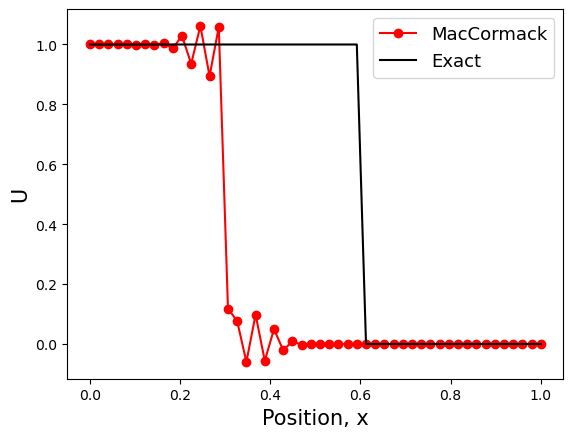}
\caption{$\mu=1000$}
\label{Fig:Ch_LB_1DLY_MC_mu1000}
\end{subfigure}
\vfill
\begin{subfigure}[b]{0.24\textwidth}
\centering
\includegraphics[width=\textwidth]{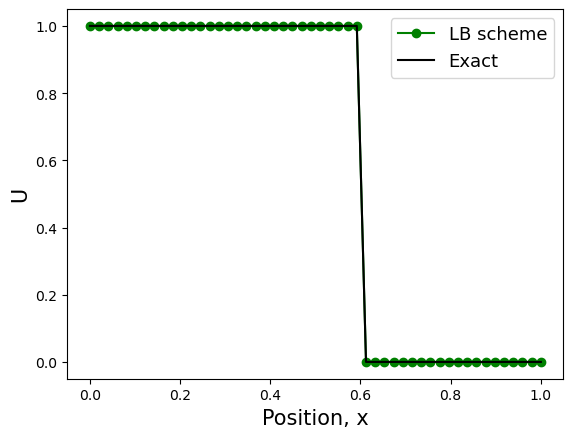}
\caption{$\mu=1$}
\label{Fig:Ch_LB_1DLY_LB_mu1}
\end{subfigure}
\hfill
\begin{subfigure}[b]{0.24\textwidth}
\centering
\includegraphics[width=\textwidth]{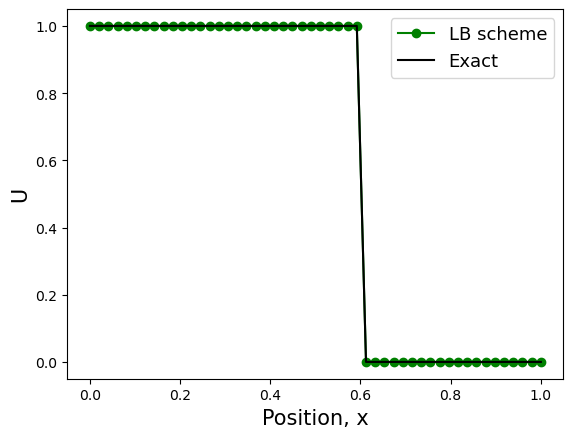}
\caption{$\mu=10$}
\label{Fig:Ch_LB_1DLY_LB_mu10}
\end{subfigure}
\hfill
\begin{subfigure}[b]{0.24\textwidth}
\centering
\includegraphics[width=\textwidth]{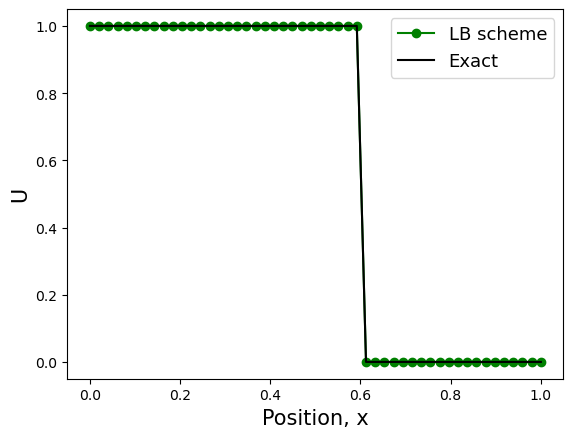}
\caption{$\mu=100$}
\label{Fig:Ch_LB_1DLY_LB_mu100}
\end{subfigure}
\hfill
\begin{subfigure}[b]{0.24\textwidth}
\centering
\includegraphics[width=\textwidth]{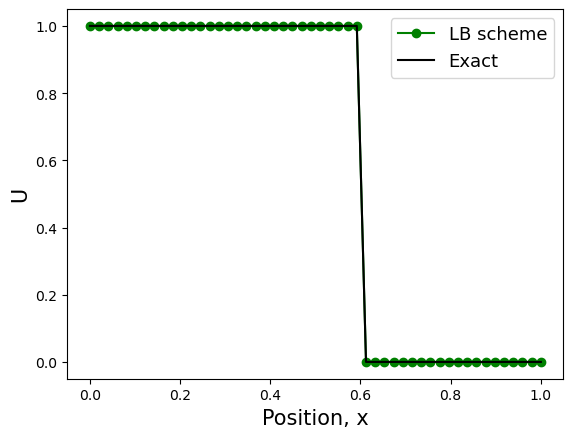}
\caption{$\mu=1000$}
\label{Fig:Ch_LB_1DLY_LB_mu1000}
\end{subfigure}
\caption{\centering Top: Extended MacCormack's method with limiter based on $U^n$(Reproduced from \cite{LeVequeYee1990}), Bottom: Our LB scheme for hyperbolic conservation laws with source terms}
\label{Fig: Ch_LB_1D LeVequeYee}
\end{figure}
 A comparison of numerical solutions reproduced from LeVeque and Yee \cite{LeVequeYee1990} and numerical solutions obtained from our LB scheme is shown in figure \ref{Fig: Ch_LB_1D LeVequeYee} at $T=0.3$ for different values of $\mu$. The MacCormack's method suffers from spurious numerical convection for $\mu$ as small as $100$, while our LB scheme is devoid of the effects of spurious numerical convection until $\mu=1000$. We observe numerical convection in LB scheme for $\mu \geq 10000$ (not shown in figure), and this validates the remark \ref{Ch_LB_source_remark} that our scheme is suitable when $S(U)=\mathcal{O}(\epsilon)$. \\
Hence, for this problem, we can infer that $\epsilon=\mathcal{O}(k\Tilde{\mu})$ where $\Tilde{\mu}$ represents the value of $\mu$ upto which the method of nullification of numerical convection works. Thus,   $\epsilon=\mathcal{O}(k10^3)$ for some constant $k<10^{-3}$.    

\subsubsection{Two dimensional discontinuity}
\label{Ch_LB_Numres_source_2D}
\begin{figure}[!h]
\centering
\begin{subfigure}[b]{0.24\textwidth}
\centering
\includegraphics[width=\textwidth]{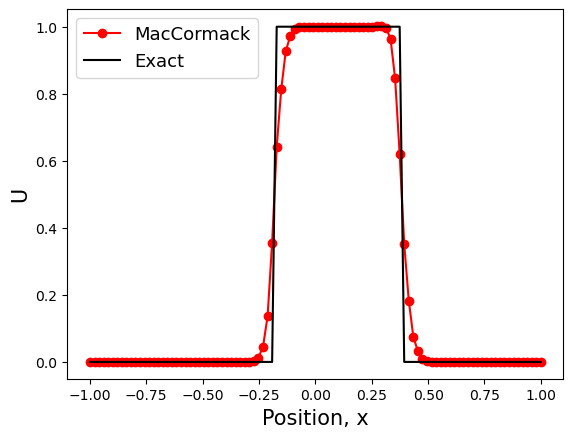}
\caption{$\mu=1$}
\label{Fig:Ch_LB_2DLY_MC_mu1}
\end{subfigure}
\hfill
\begin{subfigure}[b]{0.24\textwidth}
\centering
\includegraphics[width=\textwidth]{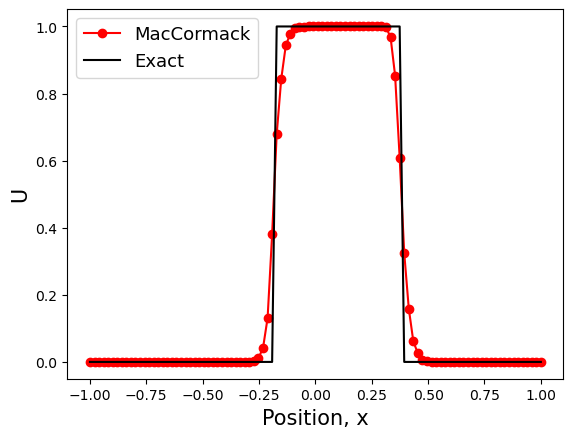}
\caption{$\mu=10$}
\label{Fig:Ch_LB_2DLY_MC_mu10}
\end{subfigure}
\hfill
\begin{subfigure}[b]{0.24\textwidth}
\centering
\includegraphics[width=\textwidth]{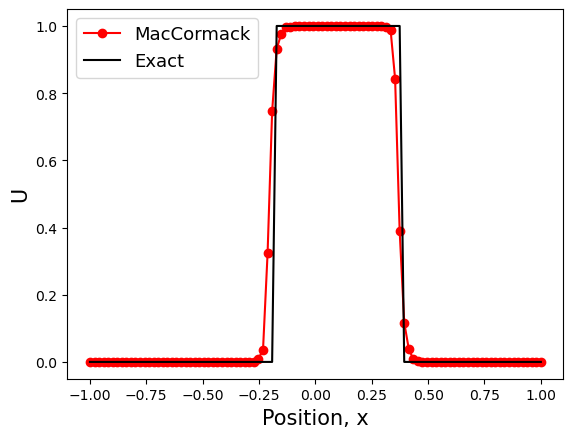}
\caption{$\mu=100$}
\label{Fig:Ch_LB_2DLY_MC_mu100}
\end{subfigure}
\hfill
\begin{subfigure}[b]{0.24\textwidth}
\centering
\includegraphics[width=\textwidth]{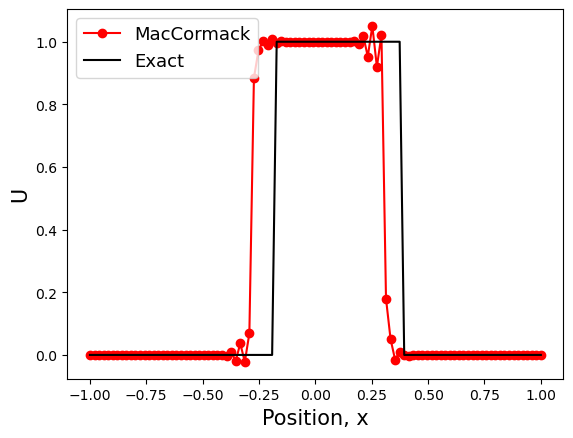}
\caption{$\mu=500$}
\label{Fig:Ch_LB_2DLY_MC_mu1000}
\end{subfigure}
\vfill
\begin{subfigure}[b]{0.24\textwidth}
\centering
\includegraphics[width=\textwidth]{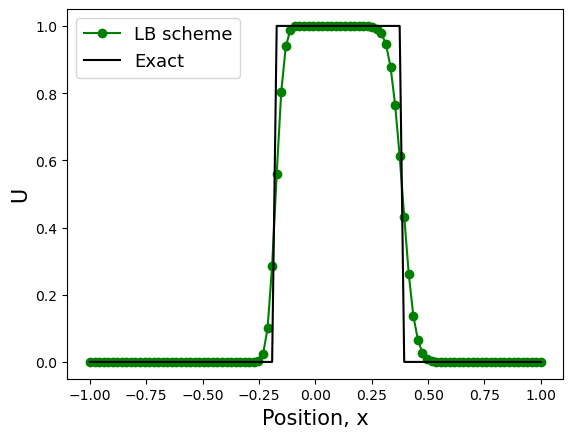}
\caption{$\mu=1$}
\label{Fig:Ch_LB_2DLY_LB_mu1}
\end{subfigure}
\hfill
\begin{subfigure}[b]{0.24\textwidth}
\centering
\includegraphics[width=\textwidth]{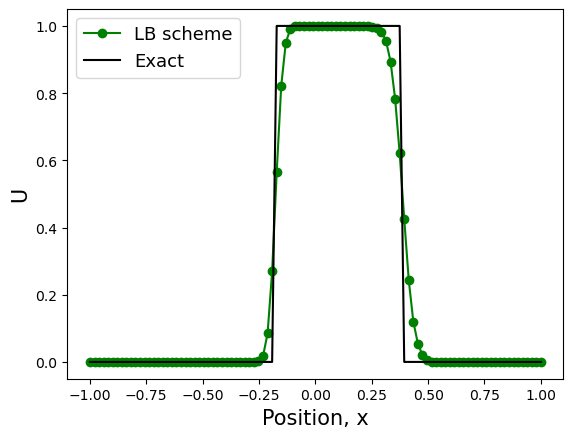}
\caption{$\mu=10$}
\label{Fig:Ch_LB_2DLY_LB_mu10}
\end{subfigure}
\hfill
\begin{subfigure}[b]{0.24\textwidth}
\centering
\includegraphics[width=\textwidth]{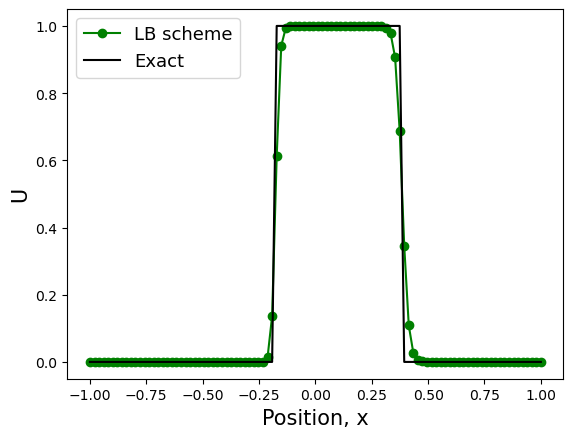}
\caption{$\mu=100$}
\label{Fig:Ch_LB_2DLY_LB_mu100}
\end{subfigure}
\hfill
\begin{subfigure}[b]{0.24\textwidth}
\centering
\includegraphics[width=\textwidth]{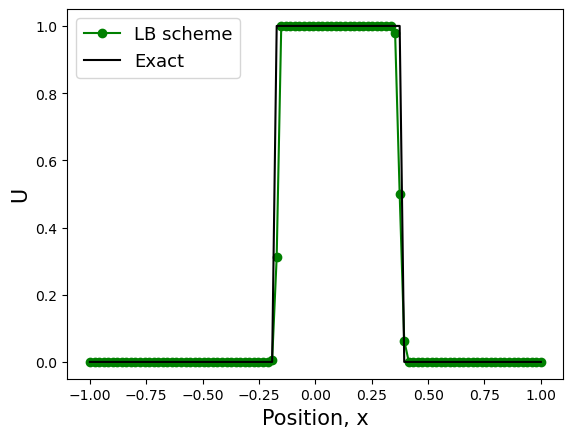}
\caption{$\mu=500$}
\label{Fig:Ch_LB_2DLY_LB_mu500}
\end{subfigure}
\caption{\centering Cross-sectional plot at $x_2=0$. \textbf{Top:} Extended MacCormack's method with limiter based on $U^n$, \textbf{Bottom:} Formulated LB scheme for hyperbolic conservation laws with source terms}
\label{Fig: Ch_LB_2D LeVequeYee}
\end{figure}

We introduce a variant of LeVeque and Yee \cite{LeVequeYee1990}'s problem in two dimensions, to understand the effect of $\epsilon$ on numerical convection. The domain is $[-1,1] \times [-1,1] \subset \mathbb{R}^2$, and is split up into $100 \times 100$ grid points. Note that $\Delta x_1 = \Delta x_2 = \Delta x$ is same as the grid spacing used in the previous one dimensional problem. For this problem, $G^1(U)=G^2(U)=U$ and $S(U)=-\mu U (U-1) (U-\frac{1}{2})$. Initial conditions are: 
\begin{equation*}
    U(x_1,x_2,0)=\left\{ \begin{matrix}
        1 & \text{for } x_1^2 + x_2^2 \leq 0.3 \\
        0 & \text{for } x_1^2 + x_2^2 > 0.3
    \end{matrix} \right.. 
\end{equation*}
An LBM with upwind D$2$Q$5$ form for $f_q^{eq}$ and $r_q$ is utilised to obtain the numerical solution. $\lambda=\frac{\Delta x}{\Delta t}$ is chosen such that the sub-characteristic condition in \eqref{Ch_LB_Upwindeqbm_subcharcondn} is satisfied. This simplifies in this case as 
\begin{equation*}
    \det \left(\begin{matrix}
        \lambda-1 & -1 \\ -1 & \lambda-1
    \end{matrix}\right) \geq 0 \implies \lambda \geq 0 \text{ and } \lambda \geq 2.
\end{equation*}
Further, the time step is chosen as $\Delta t = \frac{\Delta x}{\lambda}$. We also consider $\omega=1$ for the simulation of this problem and this ensures consistency with the governing equation irrespective of the choice of $\Delta t$ (as discussed in remark \ref{Ch_LB_Consistency_omega1}). 
A comparison of numerical solutions obtained from MacCormack's method and our LB scheme is shown in figure \ref{Fig: Ch_LB_2D LeVequeYee} at $T=0.1$ for different values of $\mu$. 
It can be seen that, for $\mu=500$, the MacCormack's method suffers from spurious numerical convection while our LB scheme does not. \\ 
In the following, we make an estimation of $\mu$ up to which our method will work according to remark \ref{Ch_LB_source_remark}. For this, we use subscripts $D=1$ and $D=2$ to compare certain variables from sections \ref{Ch_LB_Numres_source_1D} and \ref{Ch_LB_Numres_source_2D} respectively. Since $\lambda_{D=1}=1$, $\lambda_{D=2}=2$, and $\Delta x$ is the same for both one and two dimensional problems, we have $\Delta t_{D=2} = \frac{\Delta t_{D=1}}{2}$. Further, since $\omega_{D=1}=\frac{\Delta t_{D=1}}{\epsilon_{D=1}}$ and $\omega_{D=2}=\frac{\Delta t_{D=2}}{\epsilon_{D=2}}$ are both equal to $1$, we have $\epsilon_{D=2}=\frac{\epsilon_{D=1}}{2}$. Thus, since $\epsilon_{D=1}$ in section \ref{Ch_LB_Numres_source_1D} is $\mathcal{O}(k10^3)$, $\epsilon_{D=2}$ for two dimensional problem is $\mathcal{O}\left(k\frac{10^3}{2}\right)$. Hence, for this problem, our LB scheme is expected to be devoid of spurious numerical convection for $\mu$ up to $\mathcal{O}(500)$, and this is validated by the numerical results shown in figure \ref{Fig: Ch_LB_2D LeVequeYee}.

\subsubsection{Three dimensional discontinuity}
\label{Ch_LB_Numres_source_3D}
\begin{figure}[t]
\centering
\begin{subfigure}[b]{0.24\textwidth}
\centering
\includegraphics[width=\textwidth]{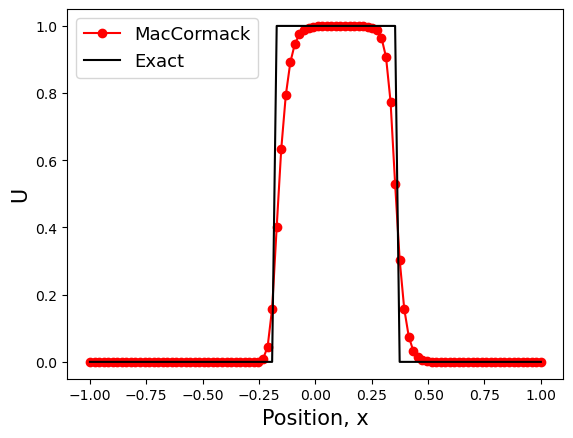}
\caption{$\mu=1$}
\label{Fig:Ch_LB_3DLY_MC_mu1}
\end{subfigure}
\hfill
\begin{subfigure}[b]{0.24\textwidth}
\centering
\includegraphics[width=\textwidth]{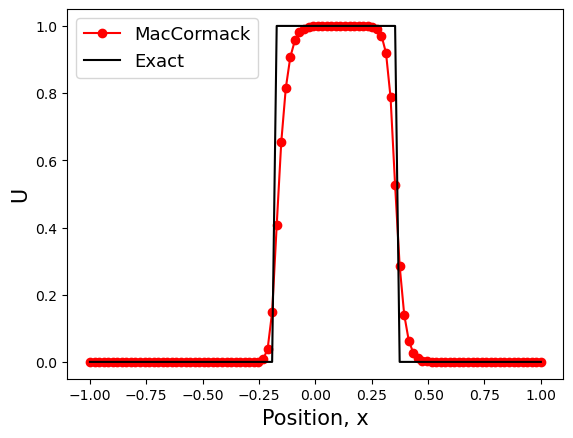}
\caption{$\mu=10$}
\label{Fig:Ch_LB_3DLY_MC_mu10}
\end{subfigure}
\hfill
\begin{subfigure}[b]{0.24\textwidth}
\centering
\includegraphics[width=\textwidth]{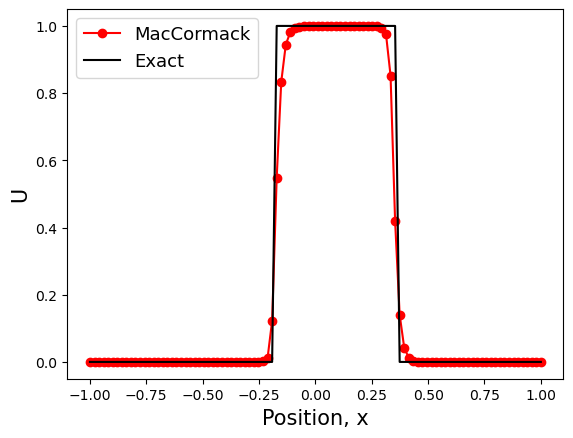}
\caption{$\mu=100$}
\label{Fig:Ch_LB_3DLY_MC_mu100}
\end{subfigure}
\hfill
\begin{subfigure}[b]{0.24\textwidth}
\centering
\includegraphics[width=\textwidth]{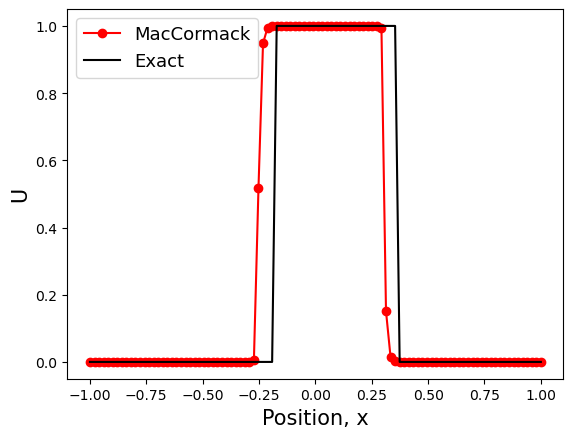}
\caption{$\mu=500$}
\label{Fig:Ch_LB_3DLY_MC_mu1000}
\end{subfigure}
\vfill
\begin{subfigure}[b]{0.24\textwidth}
\centering
\includegraphics[width=\textwidth]{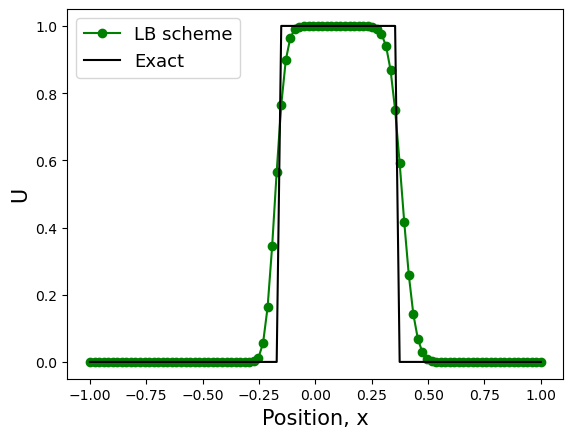}
\caption{$\mu=1$}
\label{Fig:Ch_LB_3DLY_LB_mu1}
\end{subfigure}
\hfill
\begin{subfigure}[b]{0.24\textwidth}
\centering
\includegraphics[width=\textwidth]{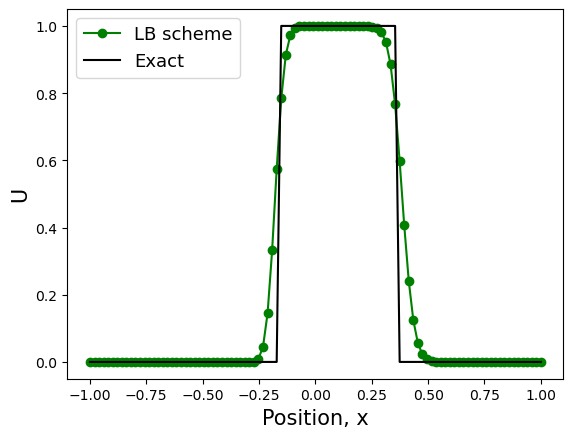}
\caption{$\mu=10$}
\label{Fig:Ch_LB_3DLY_LB_mu10}
\end{subfigure}
\hfill
\begin{subfigure}[b]{0.24\textwidth}
\centering
\includegraphics[width=\textwidth]{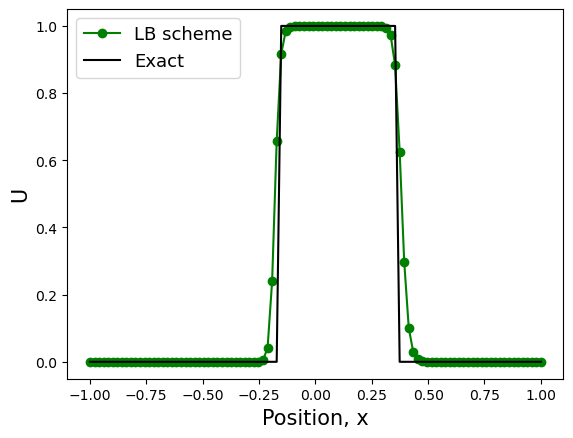}
\caption{$\mu=100$}
\label{Fig:Ch_LB_3DLY_LB_mu100}
\end{subfigure}
\hfill
\begin{subfigure}[b]{0.24\textwidth}
\centering
\includegraphics[width=\textwidth]{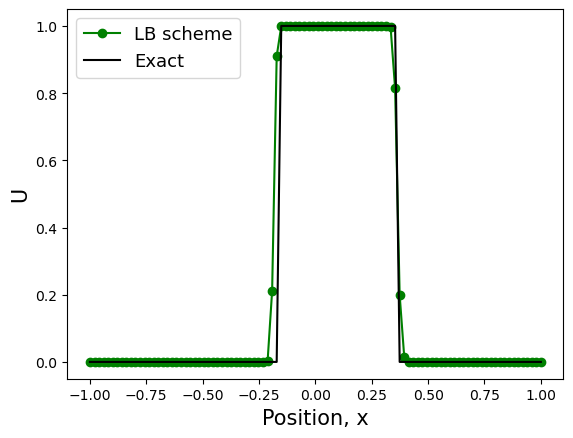}
\caption{$\mu=500$}
\label{Fig:Ch_LB_3DLY_LB_mu500}
\end{subfigure}
\caption{\centering Cross-sectional plot at $x_2,x_3=0$. \textbf{Top:} Extended MacCormack's method with limiter based on $U^n$, \textbf{Bottom:} Formulated LB scheme for hyperbolic conservation laws with source terms}
\label{Fig: Ch_LB_3D LeVequeYee}
\end{figure}

Here, we introduce a variant of LeVeque and Yee \cite{LeVequeYee1990}'s problem in three dimensions. The domain is $[-1,1] \times [-1,1] \times [-1,1] \subset \mathbb{R}^3$, and is split up into $100 \times 100 \times 100$ grid points. Note that $\Delta x_1 = \Delta x_2 = \Delta x_3 = \Delta x$ is same as the grid spacing used in the one dimensional case. For this problem, $G^1(U)=G^2(U)=G^3(U)=U$ and $S(U)=-\mu U (U-1) (U-\frac{1}{2})$. Initial conditions are: 
\begin{equation*}
    U(x_1,x_2,x_3,0)=\left\{ \begin{matrix}
        1 & \text{for } x_1^2 + x_2^2 + x_3^2 \leq 0.3 \\
        0 & \text{for } x_1^2 + x_2^2 + x_3^2> 0.3
    \end{matrix} \right.. 
\end{equation*}
An LBM with upwind D$3$Q$7$ form for $f_q^{eq}$ and $r_q$ is utilised to obtain the numerical solution. $\lambda=\frac{\Delta x}{\Delta t}$ is chosen such that the sub-characteristic condition in \eqref{Ch_LB_Upwindeqbm_subcharcondn} is satisfied. This simplifies in this case as 
\begin{equation*}
    \det \left(\begin{matrix}
        \lambda-1 & -1 & -1 \\ -1 & \lambda-1 & -1 \\ -1 & -1 & \lambda-1
    \end{matrix}\right) \geq 0 \implies \lambda \geq 0 \text{ and } \lambda \geq 3.
\end{equation*}
Further, the time step is chosen as $\Delta t = \frac{\Delta x}{\lambda}$. We also consider $\omega=1$ for the simulation of this problem and this ensures consistency with the governing equation irrespective of the choice of $\Delta t$ (as discussed in remark \ref{Ch_LB_Consistency_omega1}). A comparison of numerical solutions obtained from MacCormack's method and our LB scheme is shown in figure \ref{Fig: Ch_LB_3D LeVequeYee} at $T=0.1$ for different values of $\mu$. It can be seen that, for $\mu=500$, the MacCormack's method suffers from spurious numerical convection while our LB scheme does not. \\ 
In the following, we use subscripts $D=1$ and $D=3$ to compare certain variables from sections \ref{Ch_LB_Numres_source_1D} and \ref{Ch_LB_Numres_source_3D} respectively. Since $\lambda_{D=1}=1$, $\lambda_{D=3}=3$, and $\Delta x$ is the same for both one and three dimensional problems, we have $\Delta t_{D=3} = \frac{\Delta t_{D=1}}{3}$. Further, since $\omega_{D=1}=\frac{\Delta t_{D=1}}{\epsilon_{D=1}}$ and $\omega_{D=3}=\frac{\Delta t_{D=3}}{\epsilon_{D=3}}$ are both equal to $1$, we have $\epsilon_{D=3}=\frac{\epsilon_{D=1}}{3}$. Thus, since $\epsilon_{D=1}$ in section \ref{Ch_LB_Numres_source_1D} is $\mathcal{O}(k10^3)$, $\epsilon_{D=3}$ for three dimensional problem is $\mathcal{O}\left(k\frac{10^3}{3}\right)$. Hence, for this problem, our LB scheme is expected to be devoid of spurious numerical convection for $\mu$ up to $\mathcal{O}\left(\frac{10^3}{3}\right)$, and we observe nullification of numerical convection for $\mu$ up to $500$ in figure \ref{Fig: Ch_LB_3D LeVequeYee}.  

\subsubsection{Non-linear problem with discontinuity}
\begin{figure}[t]
\centering
\begin{subfigure}[b]{0.24\textwidth}
\centering
\includegraphics[width=\textwidth]{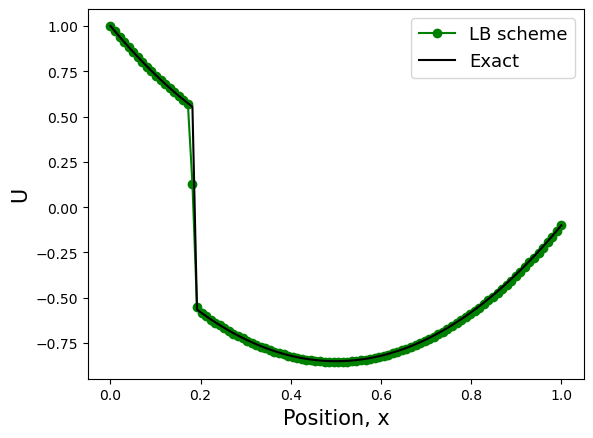}
\caption{$\mu=1$}
\label{Fig:Ch_LB_1DNL_LB_mu1}
\end{subfigure}
\hfill
\begin{subfigure}[b]{0.24\textwidth}
\centering
\includegraphics[width=\textwidth]{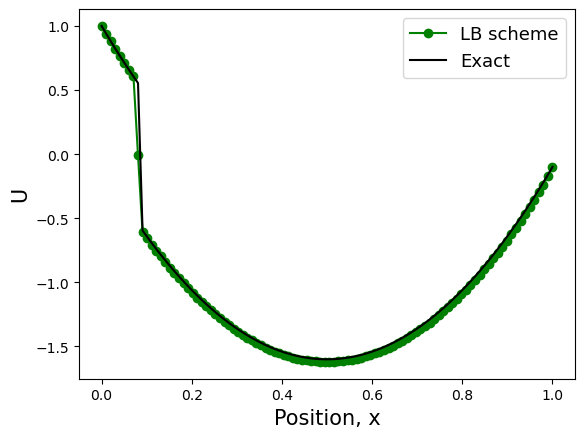}
\caption{$\mu=2$}
\label{Fig:Ch_LB_1DNL_LB_mu2}
\end{subfigure}
\hfill
\begin{subfigure}[b]{0.24\textwidth}
\centering
\includegraphics[width=\textwidth]{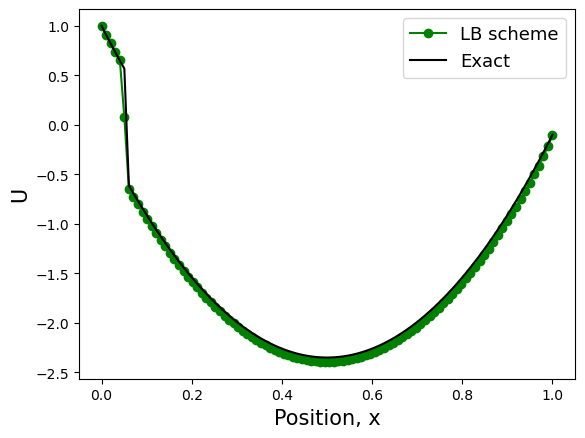}
\caption{$\mu=3$}
\label{Fig:Ch_LB_1DNL_LB_mu3}
\end{subfigure}
\hfill
\begin{subfigure}[b]{0.24\textwidth}
\centering
\includegraphics[width=\textwidth]{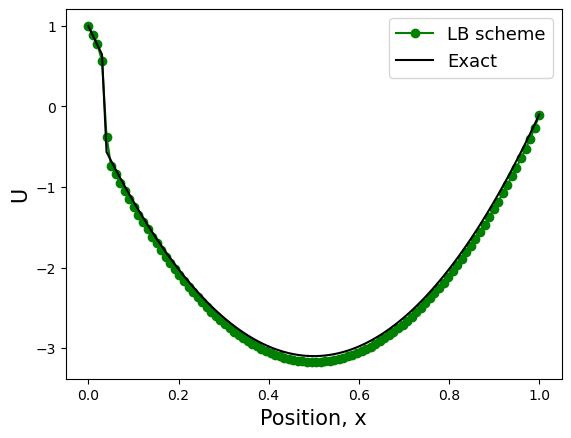}
\caption{$\mu=4$}
\label{Fig:Ch_LB_1DNL_LB_mu4}
\end{subfigure}
\vfill
\begin{subfigure}[b]{0.24\textwidth}
\centering
\includegraphics[width=\textwidth]{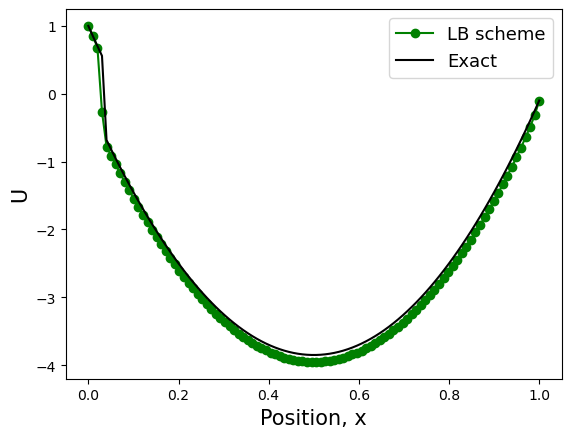}
\caption{$\mu=5$}
\label{Fig:Ch_LB_1DNL_LB_mu5}
\end{subfigure}
\hfill
\begin{subfigure}[b]{0.24\textwidth}
\centering
\includegraphics[width=\textwidth]{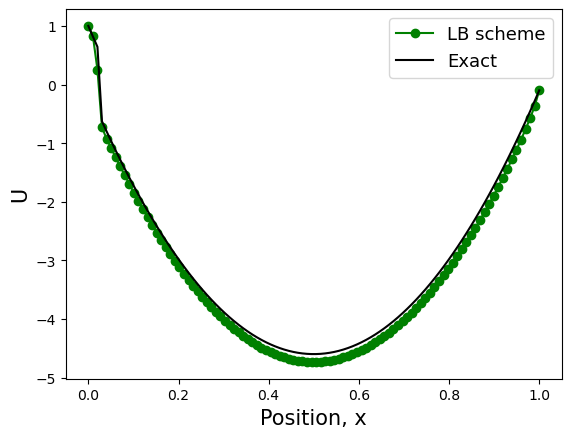}
\caption{$\mu=6$}
\label{Fig:Ch_LB_1DNL_LB_mu6}
\end{subfigure}
\hfill
\begin{subfigure}[b]{0.24\textwidth}
\centering
\includegraphics[width=\textwidth]{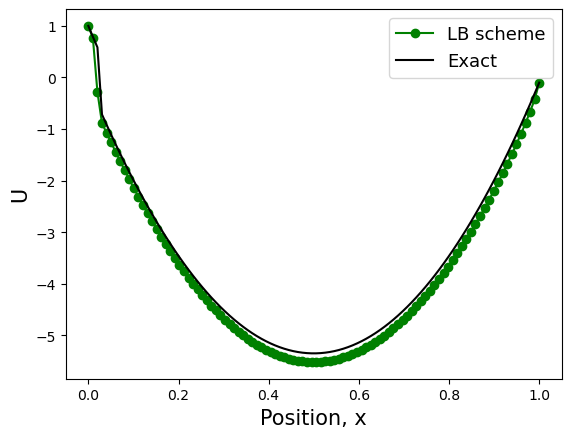}
\caption{$\mu=7$}
\label{Fig:Ch_LB_1DNL_LB_mu7}
\end{subfigure}
\hfill
\begin{subfigure}[b]{0.24\textwidth}
\centering
\includegraphics[width=\textwidth]{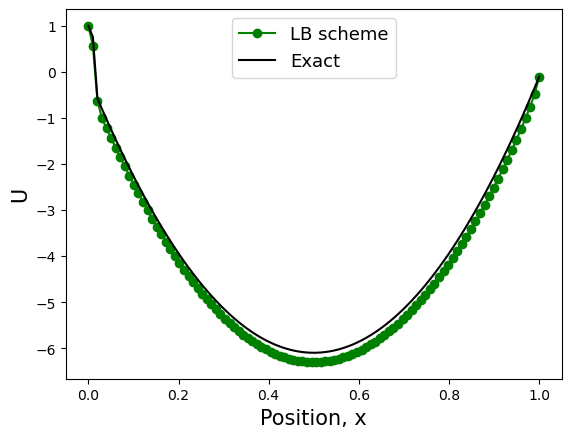}
\caption{$\mu=8$}
\label{Fig:Ch_LB_1DNL_LB_mu8}
\end{subfigure}
\caption{\centering LB scheme for non-linear problem with different values of $\mu$}
\label{Fig: Ch_LB_NLEmbid}
\end{figure}

This is a variant of the problem from Embid, Goodman and Majda \cite{EGM1984}. The domain is $[0,1] \subset \mathbb{R}$, and is split up into $100$ evenly spaced grid points. The flux function $G^1(U)=\frac{1}{2}U^2$ is non-linear and $S(U)=\mu (6x-3)U$. Boundary conditions are $U(x_1=0,t)=1$ and $U(x_1=1,t)=-0.1, \ \forall t$. For numerical simulation of this steady problem, $500$ iterations are utilised with the initialisation
\begin{equation*}
    U(x_1,0)=\left\{ \begin{matrix}
        1 & \text{for } x_1 \leq 0.1 \\
        -1 & \text{for } x_1 > 0.1
    \end{matrix} \right.. 
\end{equation*} 
 The numerical solutions obtained using LB scheme plotted against the exact solution, for different values of $\mu$, are shown in \cref{Fig: Ch_LB_NLEmbid}. It is seen that the numerical method correctly locates the discontinuities for different values of $\mu$.   

\subsection{D2Q9 model of LBM}
\begin{figure}[t]
\centering
\begin{subfigure}[b]{0.3\textwidth}
\centering
\includegraphics[width=\textwidth]{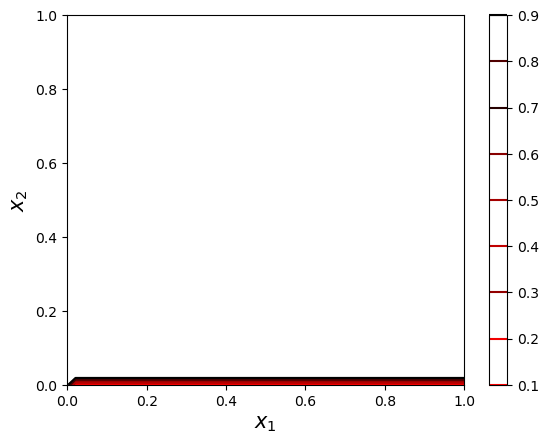}
\caption{$\theta=0$}
\label{Fig: Ch_LB_D2Q9_0}
\end{subfigure}
\hfill
\begin{subfigure}[b]{0.3\textwidth}
\centering
\includegraphics[width=\textwidth]{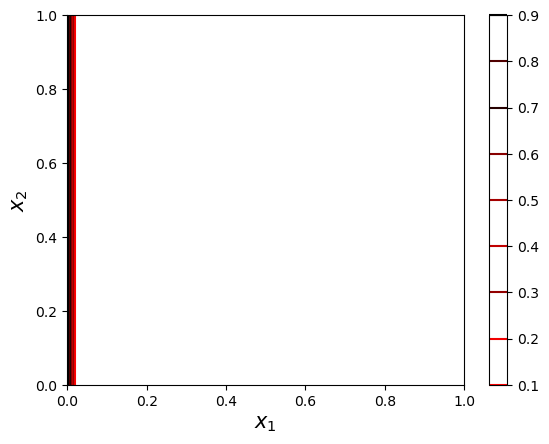}
\caption{$\theta=\frac{\pi}{2}$}
\label{Fig: Ch_LB_D2Q9_90}
\end{subfigure}
\hfill
\begin{subfigure}[b]{0.3\textwidth}
\centering
\includegraphics[width=\textwidth]{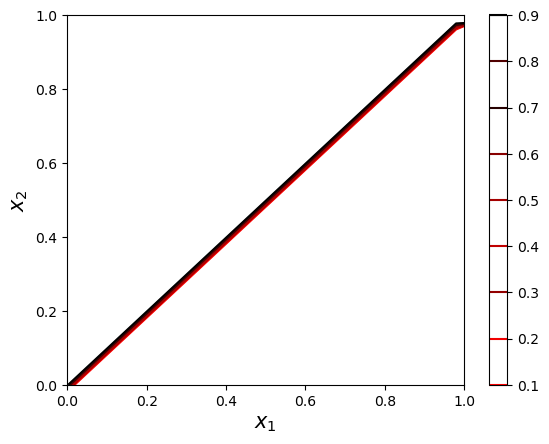}
\caption{$\theta=\frac{\pi}{4}$}
\label{Fig: Ch_LB_D2Q9_45}
\end{subfigure}
\caption{\centering Discontinuities along coordinate and diagonal-to-coordinate directions captured exactly due to upwinding} 
\label{Fig: Ch_LB_D2Q9}
\end{figure}

In this section, we show the diagonal upwinding nature of our D2Q9 model of LBM. For this, we consider a standard two-dimensional linear problem from \cite{Spekreijse1987}. The domain is $[0,1] \times [0,1] \subset \mathbb{R}^2$, and is split up into $50 \times 50$ grid points. Here $\Delta x_1 = \Delta x_2 = \Delta x$. The flux functions are $G^1(U)=aU$, $G^2(U)=bU$ where $a=\cos{\theta}$, $b=\sin{\theta}$ and $\theta \in \left(0,\frac{\pi}{2}\right)$. Boundary conditions are: 
\begin{equation*}
    \begin{matrix}
        U(0,x_2,t)=1 & \text{for } 0<x_2<1, \\
        U(x_1,0,t)=0 & \text{for } 0<x_1<1,
    \end{matrix}, \ \forall t.
\end{equation*}
Exact solution is:
\begin{equation*}
    \begin{matrix}
        U(x_1,x_2,t)=1 & \text{for } bx_1-ax_2<0, \\
        U(x_1,x_2,t)=0 & \text{for } bx_1-ax_2>0,
    \end{matrix}, \ \forall t.
\end{equation*}
It can be noted that the problem is steady. An LBM with D$2$Q$9$ equilibrium functions \eqref{Ch_LB_D2Q9_Eqbm fn} is utilised to obtain the numerical solution. For this problem, we run $1000$ iterations of our LBM before presenting the steady state solution. $\lambda$ is chosen such that the sub-characteristic condition obtained by imposing positivity of numerical diffusion coefficient in \eqref{Ch_LB_D2Q9_CE} is satisfied. Further, we consider $\omega=1$ for the simulation of this problem and this ensures consistency with the governing equation (as discussed in remark \ref{Ch_LB_Consistency_omega1}). \\
The numerical solutions for $\theta = 0 \text{ and } \theta=\frac{\pi}{2}$ obtained by choosing the fluxes $G^{\gamma}=G^{\zeta}= 0$ (thereby replicating a standard D$2$Q$5$ upwind model), are shown in figures \ref{Fig: Ch_LB_D2Q9_0} and \ref{Fig: Ch_LB_D2Q9_90} respectively. The numerical solution for $\theta = \frac{\pi}{4}$ obtained by choosing $G^{\alpha}=G^{\beta}=0$, is shown in figure \ref{Fig: Ch_LB_D2Q9_45}. It can be seen from these results that, for a specific partition of total flux between coordinate and diagonal-to-coordinate directions, the D2Q9 model captures discontinuities aligned with $x_1, x_2$ and diagonal directions exactly.

\section{Summary and conclusions}
\label{Ch_LB_SecConc}
The following are the major highlights of the paper.
\begin{itemize}
    \item An LBE is derived by semi-implicit discretisation of VKE, and its relaxation factor is compared with that of the usual LBE obtained by explicit discretisation of VKE. 
    \item The usual condition on $\hat{\omega}$ enforced by positivity of numerical diffusion coefficient in Chapman-Enskog expansion is $0<\hat{\omega}<2$. On the other hand, the properties such as H-inequality, total variation boundedness and positivity enforce the stronger constraint $0<\hat{\omega}\leq1$. By construction, the LBE that we derived by semi-implicit discretisation of VKE naturally satisfies this stronger condition. 
    \item Macroscopic finite difference form of the LBEs is derived, and it is utilised in establishing consistency of LBEs with the hyperbolic system, and in showing the total variation boundedness and the positivity of LBM. 
    \item Small numerical diffusion and better order of accuracy are realisable for $1<\hat{\omega}<2$ in the case of LBE derived by explicit discretisation of VKE. 
    \item The LBM framework is extended to hyperbolic conservation laws with source terms and the spurious numerical convection due to imbalance between convection and source terms is removed by suitable modelling of $r_q$.  The resulting method not only leads to {\em well-balancing} but also is effective for source terms of significant stiffness.  
    \item A D$2$Q$9$ model of our LBM framework allows upwinding along diagonal directions, in addition to the usual upwinding along co-ordinate directions, resulting in better multidimensional behaviour. 
\end{itemize}

\section*{CRediT author statement}
\noindent \textbf{Megala Anandan}: Conceptualization, Methodology, Formal analysis, Software, Validation, Investigation, Resources, Writing- Original draft, Reviewing and Editing, Visualization. \\
\textbf{S. V. Raghurama Rao}: Resources, Writing- Reviewing and Editing, Supervision. 

\bibliographystyle{acm}
\bibliography{references}  
\end{document}